\documentclass[12pt]{article}
\usepackage{amsmath,amsthm,amssymb,diagrams}

\newtheorem{theorem}{Theorem}[section]
\newtheorem{corollary}{Corollary}[section]
\newtheorem{lemma}{Lemma}[section]
\newtheorem{proposition}{Proposition}[section]
\theoremstyle{remark}

\theoremstyle{remark}

\theoremstyle{remark}
\newtheorem{remark}{Remark}[section]

\newcommand{\di}{\partial}
\newcommand{\D}{\mathcal D}
\newcommand{\Ffin}{\mathcal F_{\mathrm{fin}}(\mathcal D)}
\newcommand{\la}{\langle}
\newcommand{\ra}{\rangle}

\begin{document}

\makeatletter\@addtoreset{equation}{section}

\begin{center}{\Large \bf
  Meixner class of  non-commutative generalized stochastic processes with freely independent values I. A characterization
}\end{center}

{\large Marek Bo\.zejko}\\
Instytut Matematyczny, Uniwersytet Wroc{\l}awski, Pl.\ Grunwaldzki 2/4, 50-384 Wroc{\l}aw, Poland\\
e-mail: \texttt{bozejko@math.uni.wroc.pl}\vspace{2mm}

{\large Eugene Lytvynov}\\ Department of Mathematics,
Swansea University, Singleton Park, Swansea SA2 8PP, U.K.\\
e-mail: \texttt{e.lytvynov@swansea.ac.uk}\vspace{2mm}

{\small

\begin{center}
{\bf Abstract}
\end{center}

\noindent  Let $T$ be an underlying space with a non-atomic measure $\sigma$ on it (e.g.\ $T=\mathbb R^d$ and $\sigma$ is the Lebesgue measure). We introduce and study a class of non-commutative generalized stochastic processes, indexed by points of $T$, with freely independent values. Such a process (field), $\omega=\omega(t)$, $t\in T$, is given a rigorous meaning  through smearing out with  test functions on $T$, with $\int_T \sigma(dt)f(t)\omega(t)$ being a (bounded) linear operator in a full Fock space. We define a set $\mathbf{CP}$ of all continuous polynomials of $\omega$, and then define a con-commutative $L^2$-space $L^2(\tau)$ by taking the closure of $\mathbf{CP}$ in the norm
$\|P\|_{L^2(\tau)}:=\|P\Omega\|$, where $\Omega$ is the vacuum in the Fock space. Through  procedure of orthogonalization of polynomials, we construct a unitary  isomorphism between $L^2(\tau)$ and a (Fock-space-type) Hilbert space $\mathbb F=\mathbb R\oplus\bigoplus_{n=1}^\infty L^2(T^n,\gamma_n)$, with explicitly given measures $\gamma_n$. We identify the Meixner class as those processes for which the procedure of orthogonalization leaves the set $\mathbf {CP}$ invariant. (Note that, in the general case, the projection of a continuous monomial of oder $n$ onto the $n$-th chaos need not remain a continuous polynomial.) Each element of the Meixner class is characterized by two continuous functions $\lambda$ and $\eta\ge0$ on $T$, such  that, in the  $\mathbb F$ space,  $\omega$ has representation $\omega(t)=\di_t^\dag+\lambda(t)\di_t^\dag\di_t+\di_t+\eta(t)\di_t^\dag\di^2_t$, where $\di_t^\dag$ and $\di_t$ are the usual creation and annihilation operators  at point $t$.

 } \vspace{2mm}

\section{Introduction} In his classical work \cite{Meixner}, Meixner searched for all probability measures $\mu$ on $\mathbb R$ with infinite support whose system of monic orthogonal polynomials $(p^{(n)})_{n=0}^\infty$ has an (exponential) generating function of the exponential type:
\begin{equation}\label{ufytfdtk}\sum_{n=0}^\infty \frac{p^{(n)}(t)}{n!}\,z^n=\exp(t\psi(z)+\phi(z))=\sum_{k=0}^\infty\frac1{k!}\, (t\psi(z)+\phi(z))^k.\end{equation} Meixner discovered that this (essentially) holds if and only if there exist $\lambda\in\mathbb R$ and $\eta\ge0$ such that the polynomials $(p^{(n)})_{n=0}^\infty$ satisfy the recursive relation
\begin{equation}\label{dtdstj} tp^{(n)}(t)=p^{(n+1)}(t)+\lambda np^{(n)}(t)+
(n+\eta n(n-1))p^{(n-1)}(t).\end{equation}
(We refer to  \cite{S} for a modern presentation of this result.) From \eqref{dtdstj}
one concludes that the measure $\mu$ can be either Gaussian, or Poisson, or gamma, or Pascal (negative binomial), or Meixner.
We may now introduce  in $L^2(\mathbb R,\mu)$ creation (raising) and annihilation (lowering) operators through
$\di^\dag p^{(n)}:=p^{(n+1)}$ and
$\di p^{(n)}:=np^{(n-1)}$, respectively.
Then, by \eqref{dtdstj}, the action of the operator of multiplication by $t$ in $L^2(\mathbb R,\mu)$ has a representation
\begin{equation}\label{txsjusxt} t\cdot=\di^\dag+\lambda \di^\dag\di +\di+\eta\di^\dag\di\di.\end{equation}

 Since Meixner's laws are infinitely divisible, they appear as distributions of increments of corresponding L\'evy processes.
 These   are exactly Brownian motion,   Poisson, gamma,  Pascal, and Meixner processes.
 Note the first two of these  processes correspond to the case $\eta=0$, while the latter three correspond to $\eta>0$.  We will refer to all of them as  the Meixner class of L\'evy processes. From numerous applications of these processes let us mention that, for $\eta>0$, they naturally appear  in the study of  a realization of the renormalized  square of white noise, see \cite{AFS,Ly3,sn} and the references therein.

In \cite{Ly1} (see also \cite{bm,KL,KSSU,Rod}), Meixner-type generalized stochastic processes with independent values were constructed and studied.  More precisely, consider a standard triple of the form
$ S\subset L^2(\mathbb R,dt)\subset S'$, where
$S$ is a nuclear space of smooth functions, and $S'$ is the dual of $S$ with respect to the central space $L^2(\mathbb R,dt)$, i.e., $S'$ is a space of generalized functions (distributions). Let $\lambda,\eta$ be parameters as in \eqref{dtdstj}, or even more generally, let $\lambda(\cdot)$  and $\eta(\cdot)$ be smooth functions on $\mathbb  R$, which give, at each $t\in\mathbb  R$, parameters $\lambda(t),\eta(t)$. Then, there exists a probability measure $\mu$ on the space $S'$ which is a generalized  stochastic process with independent values (in the sense of \cite{GV}), and the operator of multiplication by a monomial $\la f,\omega\ra$, $f\in S$, $\omega\in S'$, has a representation
$$ \la f,\omega\ra\cdot=\int_{\mathbb R}dt\,f(t)(\di_t^\dag+\lambda(t)\di_t^\dag\di_t+\di_t+\eta(t)\di_t^\dag\di_t\di_t),$$
that is,
\begin{equation}\label{ctyd} \omega(t)\cdot=\di_t^\dag+\lambda(t)\di_t^\dag\di_t+\di_t+\eta(t)\di_t^\dag\di_t\di_t \end{equation}
(compare with \eqref{txsjusxt}). In \eqref{ctyd}, the operators $\di_t^\dag$ and $\di_t$ are defined by analogy with the one-dimensional case, although on infinite-dimensional orthogonal polynomials of $\omega\in S'$, so that $\di_t^\dag$ and $\di_t$ are the usual creation and annihilation operators  at point $t$.  As a result, one has a unitary isomorphism between the $L^2$-space $L^2(S',\mu)$ and some Hilbert space $ F=\bigoplus_{n=0}^\infty F^{(n)}$, where $F^{(0)}=\mathbb R$, while for each $n\in\mathbb N$,
$  F^{(n)}=L_{\mathrm{sym}}^2(\mathbb R^n,\theta_{n})$ --- the space of all symmetric functions on $\mathbb R^n$ which are square integrable with respect to some measure $\theta_{n}$ (depending on $\lambda$ and $\eta$). In the special case where $\eta\equiv0$, the space $ F$ reduces to the usual symmetric Fock space over $L^2(\mathbb R,dt)$, whereas in the general case the space $ F$ is wider than the Fock space, which is why, in \cite{Ly1},  $ F$ was called an extended Fock space.

As follows from \cite{blm,Ly2}, the Meixner class  may be characterized between all generalized stochastic processes with independent values as exactly those processes whose
orthogonal polynomials remain {\it continuous} polynomials. Recall that, in infinite dimensions, orthogonalization of polynomials  means: first,
 decomposing the $L^2$-space into the infinite orthogonal sum of its subspaces generated by polynomials, and   second, taking the projection of each monomial of order $n$ onto the $n$-th space. This is why, although the initial monomials are continuous functions of $\omega\in S'$, their orthogonal projections do not need to retain this property.  The result of
\cite{blm,Ly2} also means that it is only for the Meixner-type processes that the multiplication operator $\omega(t)\cdot$ can be represented
through the operators $\di_t^\dag$, $\di_t$.

In free probability, Meixner's systems of polynomials (on $\mathbb R$) were introduced by Anshelevich  \cite{a1} and Saitoh, Yoshida \cite{SY}. (In fact, such polynomials had already occurred in many places in the literature even before \cite{a1,SY}, see \cite[p.~62]{BB} and \cite[p.~864]{a2} for  bibliographical references.) The free Meixner polynomials $(q^{(n)})_{n=0}^\infty$ have a  (usual) generating function of the resolvent type:
 \begin{equation} \label{ygfufydts}\sum_{n=0}^\infty q^{(n)}(t)z^n=(1-(t\psi(z)+\phi(z)))^{-1}=\sum_{k=0}^\infty (t\psi(z)+\phi(z))^k\end{equation}
(compare with \eqref{ufytfdtk}).
Recall the following notation from $q$-analysis: for each $q\in[-1,1]$, we define $[0]_q:=0$ and $[n]_q:=1+q+q^2+\dots+q^{n-1}$ for $n\in\mathbb N$. In particular, for $q=0$, we have $[0]_0=0$ and $[n]_0=1$ for $n\in\mathbb N$. Then, by \cite{a1}, equality \eqref{ygfufydts} (essentially) holds if and only if there exist $\lambda\in\mathbb R$ and $\eta\ge0$ such that the polynomials $(q^{(n)})_{n=0}^\infty$ satisfy the recursive relation
\begin{equation}\label{yted6}tq^{(n)}(t)=q^{(n+1)}(t)+\lambda [n]_0q^{(n)}(t)+
([n]_0+\eta [n]_0[n-1]_0)q^{(n-1)}(t),\end{equation}
or, equivalently, equality \eqref{txsjusxt} holds in which $\di^\dag$ and $\di$ are defined through
$\di^\dag q^{(n)}:=q^{(n+1)}$ and $\di q^{(n)}:=[n]_0q^{(n-1)}$.

Each measure of orthogonality of a free Meixner system of polynomials (which has an infinite support) is freely infinitely divisible,
and therefore there exists a corresponding free L\'evy processes. A characterization of these processes in terms of a regression problem was given in \cite{BB}. These processes also   appeared
 in the study of a realization of the renormalized  square of free white noise \cite{sn}.
A deep study of free Meixner polynomials of $d$ ($d\in\mathbb N$) non-commutative variables has been carried out by Anshelevich in
 \cite{a3,a2,a5,a7}.

The aim of the present paper is to introduce and study  the Meixner class of  non-commutative generalized stochastic processes with freely independent values, or equivalently Meixner-type free polynomials of infinitely-many (non-commutative) variables. We ``translate'' the aforementioned results of the theory  of classical generalized stochastic processes with independent values into the language of free probability.
In particular, we derive  representation \eqref{ctyd}  for these processes in which $\di^\dag_t$ and $\di_t$ are the creation and annihilation operators, as in the full Fock space, at point $t$.
The main result of the paper---Theorem \ref{tyfxdzxsqx}---is the characterization of the Meixner class as  exactly those non-commutative generalized  stochastic processes with freely independent values whose orthogonal polynomials are continuous in $\omega$.

It should be stressed that, generally speaking, the orthogonal polynomials we consider resemble one-dimensional free Meixner polynomials only in the infinitesimal sense, i.e.,  at each point of the underlying space.

The paper is organized as follows. We start, in Section~\ref{cfrd}, with a discussion of processes of Gauss--Poisson type. We fix an underlying space $T$ and a non-atomic measure $\sigma$ on it. (Although the most importanant case is when $T$ is either  $\mathbb R$ or $[0,\infty)$ and $\sigma$ is the Lebesgue measure, we prefer to deal with a general space to stress that its structure does not play any significant role.)
We fix a function $\lambda\in C(T)$, and consider a process (noise) of the form $\omega(t)=\di_t^\dag+\di_t+\lambda(t)\di_t^\dag\di_t$ in the full Fock space over $L^2(T,\sigma)$.
A sense to this process is  given through smearing out with a test function $f$ on $T$.
We introduce a free expectation $\tau$ and the corresponding (non-commutative) $L^2$-space $L^2(\tau)$. In terms of the expansion through orthogonal polynomials of $\omega$, the space $L^2(\tau)$ is unitarily isomorphic to the original Fock space. We prove that the procedure of orthogonalization in $L^2(\tau)$ is equivalent to the procedure of free Wick (normal)  ordering of the operators $\di^\dag_t$ and  $\di_t$. This, in particular, generalizes a corresponding result of \cite[p.~137]{bks}, which was proved in the Gaussian case, i.e.,  when $\lambda\equiv0$ (compare also with  \cite[p.~186]{a4}). We note, however, that in \cite{bks}, the authors did not use the Wick ordering in the infinitesimal sense, which is only possible when   $\lambda\equiv0$. We then derive  theorems giving a Wick rule for the product $\omega(t_1)\dotsm\omega(t_n)$, as well as a Wick rule for a product of Wick products. The latter theorems present a free counterpart of results of \cite{LRS}, see also \cite[Proposition~6]{a4} for a $q$-case.

In Section~\ref{serts}, we study (quite) general non-commutative generalized  stochastic processes with freely independent values.
They are described by assigning to each $t\in T$, a compactly supported probability measure $\mu(t,ds)$ on $\mathbb R$, so that $\mu(t,\{0\})$
 is the diffusion coefficient of the process, while outside zero $\nu(t,ds):=\frac1{s^2}\mu(t,ds)$  is the L\'evy measure of ``jumps'' at point $t$ (compare with \cite{bnt1,bnt2,bnt3}).  We prove that the set of continuous polynomials of $\omega$ is dense in the corresponding space $L^2(\tau)$, introduce orthogonal polynomials, decompose any element of
 $L^2(\tau)$ into  an infinite sum of orthogonal polynomials, and thus derive a unitary isomorphism between $L^2(\tau)$ and an extended full Fock space $\mathbb F=\bigoplus_{n=0}^\infty\mathbb F^{(n)}$, where $\mathbb F^{(n)}=L^2(T^n,\gamma_{n})$ with some measure $\gamma_{n}$ on $T^n$. We also present an explicit form of the action of the operators of (left) multiplication by $\la f,\omega\ra$ realized in the space $\mathbb F$.  These operators have a clear Jacobi-field structure (compare with \cite{ber1,blm,bruning,Ly}).
To derive our results, we produce an expansion of $L^2(\tau)$ in multiple stochastic integrals, by analogy with the Nualart and Schoutens result \cite{NS} in the classical case. In fact, Anshelevich   \cite{a4} extended the result of \cite{NS} to the case of general $q$-L\'evy processes. Comparing our result in this section with that of \cite{a4}, we note that, first, we do not assume the process to be stationary, i.e., we allow the L\'evy measure to depend on $t$, and second, what is much more important, our main results in this
section---Theorems~\ref{drse} and
\ref{yufytrkjh}---are new even in the stationary case (when $q=0$).

Finally, in Section~\ref{gyder5asr}, we derive the Meixner class  of free processes as exactly those non-commutative generalized  stochastic processes with freely independent values for which orthogonal polynomials are continuous in $\omega$,  and thus we derive a counterpart of formula  \eqref{ctyd} in the free case.

In the second part of this paper, which is currently in preparation, we will discuss the generating function for the orthogonal polynomials of $\omega$ from the Meixner class and  other related problems, and we will also mention some open problems.

\section{Free Gauss-Poisson process}\label{cfrd}

Let $T$ be a locally compact, second countable Hausdorff topological space. Recall that such a space is known to be Polish. A subset  of $T$ is called bounded if it is relatively compact in $T$.
We will additionally assume that $T$ does not possess isolated points, i.e., for every $t\in T$, there exists a sequence $\{t_n\}_{n=1}^\infty\subset T$ such that $t_n\ne t$ for all $n\in\mathbb N$, and $t_n\to t$ as $n\to\infty$.
We denote by $\mathcal  B(T)$  the Borel $\sigma$-algebra in $T$, and by  ${\cal B}_0(T)$ the collection of all relatively compact sets from ${\cal B}(T)$.
Let  $\D:=C_0(T)$ denote the set  of all real-valued continuous
functions on $T$ with compact support. Analogously, we define $\D^{(n)}:=C_0(T^n)$, $n\in \mathbb{N}$, and ${\D}^{(0)}:=\mathbb{R}$.

For a real, separable Hilbert space $\cal H$ we denote by  $\mathcal F({\mathcal H})$ the full Fock space over $\cal H$, i.e., $
{\cal F({\cal H})}:=\bigoplus_{n=0}^\infty {\cal H}^{\otimes n}
$, where ${\cal H}^{\otimes0}:=\mathbb{R}$. As usual, we will identify each ${\cal H}^{\otimes n}$ with the corresponding subspace of ${\cal F({\cal H})}$. We denote by $\Ffin$ the subset of ${\cal F({\cal H})}$ consisting of all sequences $f=(f^{(0)},f^{(1)},\dots,f^{(n)},0,0,\dots)$ such that
$f^{(i)}\in{\cal D}^{(i)}$, $i=0,1,\dots,n$, $n\in \mathbb{N}_0:=\mathbb{N}\cup\{0\}$ .
The element $\Omega{:=}(1,0,0,\dots)\in \Ffin$ is called the vacuum.

 Let $\sigma$ be a Radon, non-atomic measure on
$(T,{\cal B}(T))$.
We will assume that the measure $\sigma$ satisfies $\sigma(O)>0$ for each open, non-empty set $O$ in $T$.
 Let ${\cal
H}{:=}L^2(T,\sigma)$ be the real $L^2$-space over $T$ with respect to the measure $\sigma$, and thus we get the Fock space ${\cal F}({\cal H})={\cal F}(L^2(T,\sigma))$.

For each $f\in \D$, we denote by $a^+(f)$, $a^-(f)$, and $a^0(f)$ the corresponding creation, annihilation, and neutral operators, respectively. These are bounded linear operators on $\cal F({\cal H})$ given through \begin{align*}&a^+(f)=f\otimes g^{(n)}, \quad g^{(n)}\in{\cal
H}^{\otimes n},\ n\in\mathbb N_0, \\
& a^-(f)\,g_1\otimes\cdots\otimes g_n=(f,g_1)_{\cal H}\,g_2\otimes\cdots\otimes g_n,\\
&a^0(f)\,g_1\otimes\cdots\otimes g_n=(fg_1) \otimes g_2\otimes\cdots\otimes g_n,\quad
 g_1,\dots, g_n\in {\cal H},\ n\in \mathbb{N},\\
&a^-(f)\,\Omega=a^0(f)\,\Omega=0.
\end{align*}
The operator  $a^+(f)$ is the adjoint of $a^-(f)$, whereas   $a^0(f)$ is self-adjoint. Note that $a^+(f)$ and $a^-(g)$, $f$, $g\in \D$, satisfy the free commutation relation
\begin{equation}\label{jsfhd}a^-(g)a^+(f)=(g,f)_{\cal H},\end{equation}
where, as usual, a constant  is understood as the constant times the identity operator $\mathbf 1$.

Throughout the paper, we will heavily use the following standard notations. 
For each $t\in T$, we define $\di_t$ as the annihilation operator at point $t$.
More precisely, we set $\di_t\Omega:=0$, and
for each  $f^{(n)}\in \D^{(n)}$, $n\in\mathbb N$, we set $$(\di_tf^{(n)})(t_1,\dots,t_{n-1}):=f^{(n)}(t,t_1,\dots,t_{n-1}).$$ Clearly, $\di_t f^{(n)}\in \D^{(n-1)}$. Extending by linearity, we see that $\di_t$ maps $\Ffin$ into itself.
If we introduce the ``delta-function'' $\delta_t$: $\langle\delta_t,f\rangle=f(t)$ for $f\in \D$, then the operator $\di_t$ can be thought of as $a^-(\delta_t)$.

Next, we heuristically define $\di_t^\dag$ as the creation operator at point $t$, i.e., $\di_t^\dag$ is the ``adjoint'' of $\di_t$, so that $\di_t^\dag=a^+(\delta_t)$. A rigorous meaning to formulas involving  $\di_t^\dag$ will be given through smearing with test functions. In particular, for each $f\in\mathcal D$, we get: 
\begin{equation}a^+(f)=\int_T\sigma(dt)f(t)\di_t^\dag,\quad
a^-(f)=\int_T\sigma(dt)f(t)\di_t, \quad
a^0(f)=\int_T\sigma(dt)f(t)\di_t^\dag\di_t.\label{sgdhjk}
 \end{equation}
 Note that the relation \eqref{jsfhd} can now be written down in the form \begin{equation}\label{dhol}\di_s\di_t^\dag=\delta(s,t),\end{equation} where
 \begin{equation}\label{hgycdh}\int_T\sigma(ds)\int_T\sigma(dt)\delta(s,t)f^{(2)}(s,t):=\int_T\sigma(dt)\,f^{(2)}(t,t),
 \quad f^{(2)}\in {\cal D}^{(2)}.
 \end{equation}

 We now fix $\lambda\in C(T)$ --- the space of all continuous functions on $T$, and define, for each $f\in  {\cal D}$ a self-adjoint operator $$x(f):=a^+(f)+a^-(f)+a^0(\lambda f),$$ so that  \begin{equation}\label{slgn}x(f)=\int_T\sigma(dt)\,f(t)\big(\di_t^\dag+\di_t+\lambda(t)\di_t^\dag\di_t\big).\end{equation}
 As we will see below, if $\lambda\equiv0$, then
$(x(f))_{f\in {\cal D}}$ is a free Gaussian process, and if $\lambda\equiv 1$, then
 $(x(f))_{f\in {\cal D}}$ is a free (centered) Poisson process. In view of \eqref{slgn}, we denote $$\omega(t):=\di_t^\dag+\di_t+\lambda(t)\di_t^\dag\di_t,$$ so that \eqref{slgn} becomes $$x(f)=\int_T\sigma(dt)f(t)\omega(t).$$ Thus, $\omega:=(\omega(t))_{t\in T}$ can be interpreted as the corresponding free  noise.

\begin{lemma}\label{t6yh} The vacuum vector $\Omega$ is cyclic for the operator family $(x(f))_{f\in {\cal D}}$, i.e.,
 \begin{equation*}
 \operatorname{c.l.s.}\{\Omega,\ x(f_1)\dotsm x(f_n)\Omega \mid f_1,\dots,f_n\in\mathcal D,\ n\in\mathbb N\}
 =\mathcal F(\mathcal H).
 \end{equation*}
 Here and below, $\operatorname{c.l.s.}$ stands for the closed linear span.
  \end{lemma}

\begin{proof} The statement follows by induction from the fact that we have a Jacobi field, i.e., each operator $a(f)$ has a three-diagonal structure, with $a^+(f)$, $f\in\D$, being the usual creation operators (compare with e.g.\ \cite{ber1,Ly}).\end{proof}

 We can  naturally extend  the definition of $x(f)$ to the case where $f\in B_0(T)$ --- the space of all real-valued bounded measurable functions on $T$ with compact support. Let $\mathbf A$ denote the real algebra generated by  $(x(f))_{f\in B_0(T)}$. We define a free expectation on  $\mathbf A$ by
 $$\tau(a):=\big(a\Omega, \Omega\big)_{{\cal F({\cal H})}} ,\quad a\in {\mathbf A}.$$

Recall that a set partition $\pi$ of a set $X$ is a collection of disjoint subsets of $X$ whose union equals $X$. Let $ {NC}(n)$ denote the collection of all non-crossing partitions of $\{1,\dots,n\}$, i.e., all set partitions $\pi=\{A_1,\dots,A_k\}$, $k\ge1$, of $\{1,\dots,n\}$ such that there do not exist $A_i,A_j\in\pi$, $A_i\ne A_j$, for which the following inequalities hold: $x_1<y_1<x_2<y_2$
for some $x_1,x_2\in A_i$ and $y_1,y_2\in A_j$.

For each $n\in\mathbb N$, we define a free cumulant $C^{(n)}$ as the $n$-linear mapping $C^{(n)}:B_0(T)^n\to\mathbb R$
defined recurrently by the following formula, which connects the free cumulants with moments:
\begin{equation} \label{tdtrs}\tau(x(f_1)x(f_2)\dotsm x(f_n))=\sum_{\pi\in{NC}(n)}\prod_{A\in\pi}C(A,f_1,\dots,f_n),\end{equation}
where for each $A=\{a_1,\dots,a_k\}\subset \{1,2,\dots,n\}$, $a_1<a_2<\dots<a_k$,
$$ C(A,f_1,\dots,f_n):=C^{(k)}(f_{a_1},\dots,f_{a_k}).
$$
As easily seen, $C^{(1)}\equiv0$ and
\begin{equation}\label{dthufg}
 C^{(n)}(f_1,\dots,f_n)=\int_T f_1(t)\dotsm f_n(t)\lambda^{n-2}(t)\,\sigma(dt),\quad f_1,\dots, f_n\in B_0(T),\ n\ge2.
 \end{equation}
 By \eqref{tdtrs} and \eqref{dthufg}, the expectation $\tau$ on $\mathbf A$ is tracial, i.e., for any $a,b\in\mathbf A$, $\tau(ab)=\tau(ba)$.

\begin{proposition}\label{dreshg}
Let $f_1,\dots,f_n\in B_0(T)$ be such that
\begin{equation}\label{ytfd}
f_if_j=0\quad \text{$\sigma$-a.e.\ for all\/ $1\le i<j\le n.$}\end{equation}
Then $x(f_i)$, $i=1,\dots,n$, are freely independent with respect to $\tau$. \end{proposition}

\begin{proof} 
By \eqref{dthufg} and \eqref{ytfd}, for each $k\ge2$ and any indices $i_1,\dots,i_k\in\{1,\dots,n\}$ such that $i_l\ne i_m$ for some $l,m\in\{1,\dots,k\}$, $C^{(k)}(f_{i_1},\dots,f_{i_k})=0$. Using e.g.\ \cite{Speicher}, we conclude from here the statement.
\end{proof}

Let $B_0(T)_{\mathbb C}$ denote the complexification of $B_0(T)$.
We extend $C^{(n)}$
by linearity to the $n$-linear mapping $C^{(n)}:B_0(T)_{\mathbb C}^n\to\mathbb C$. For each $f\in
B_0(T)_{\mathbb C}$, we denote $C^{(n)}(f):=C^{(n)}(f,\dots,f)$, and define the free cumulant transform
$ C(f):=\sum_{n=1}^\infty C^{(n)}(f)$,
provided that the latter series converges absolutely. By \eqref{dthufg} and the dominated convergence theorem, we get:

\begin{proposition}\label{yuftdfrl}
Let $f\in B_0(T)_{\mathbb C}$ be such that there exists $\varepsilon\in(0,1)$ for which
\begin{equation}\label{dyrdry} |f(t)|<\frac{1-\varepsilon}{\lambda(t)} \quad\text{for all }t\in T\end{equation}
(where  $\frac{1-\varepsilon}0:=+\infty$). Then
\begin{equation}\label{afyt}
C(f)=\int_T\frac{f^2(t)}{1-\lambda(t)f(t)}\,\sigma(dt).\end{equation}
\end{proposition}

\begin{remark}\label{iuoyhh}
Note that, for $f\in \D$, condition \eqref{dyrdry} is equivalent to $|f(t)|<1/|\lambda(t)|$ for all $t\in T$.
\end{remark}

For each $\Delta\in \mathcal B_0(T)$, we define
$$ x(\Delta):=x(\chi_\Delta)=\int_\Delta \sigma(dt)\omega(t),$$
where $\chi_\Delta$ denotes the indicator function of $\Delta$. Then, by Proposition~\ref{dreshg},
for any mutually disjoint sets $\Delta_1,\dots,\Delta_n\in\mathcal B_0(T)$, the operators $x(\Delta_1),\dots,x(\Delta_n)$
are freely independent, and so by analogy with the classical case (see e.g.\ \cite{GV}), we can interpret $\omega$
as a non-commutative generalized stochastic process with freely independent values.


 For each $f^{(n)}\in B_0(T^{n})$, we define a monomial of $\omega$ by  \begin{align}\langle f^{(n)}, \omega^{\otimes n}\rangle:&=\int_{T^n}\sigma(dt_1)\dotsm\sigma(dt_n)\,f^{(n)}(t_1,\dots, t_n)\omega(t_1)\dotsm\omega(t_n)\notag\\
 &=\int_{T^n}\sigma(dt_1)\dotsm\sigma(dt_n)\,f^{(n)}(t_1,\dots, t_n) (\di_{t_1}^\dag+\di_{t_1}+\lambda(t_1)\di_{t_1}^\dag\di_{t_1})\notag\\
&\qquad\times\dots\times 
 (\di_{t_n}^\dag+\di_{t_n}+\lambda(t_n)\di_{t_n}^\dag\di_{t_n}).\label{owij}\end{align} 
 In fact,
the presence of $\di_{t_i}^\dag$ in \eqref{owij} just means  the creation of a function in the $t_i$-variable, the presence of $\lambda(t_i)\di_{t_i}^\dag\di_{t_i}$ means the identification of the $t_i$-variable with the previous $t_{i-1}$-variable and additional multiplication by $\lambda(t_i)$, whereas the presence of $\di_{t_i}$ means integration in the $t_i$-variable. For example, for $f^{(4)}\in  B_0(T^4)$ and $g^{(2)}\in \D^{(2)}$,
\begin{align*}
&\bigg(\int_{T^4}\sigma(dt_1)\dotsm\sigma(dt_4)f^{(4)}(t_1,\dots,t_4)\di^\dag_{t_1}\di_{t_2}\lambda(t_3)\di_{t_3}^\dag\di_{t_3}\di_{t_4}^\dag g^{(2)}\bigg)(s_1,s_2,s_3)
\\&\qquad =\int_T\sigma(dt)\lambda(t)f^{(4)}(s_1,t,t,t)g^{(2)}(s_2,s_3).
\end{align*}

Using the Cauchy--Schwarz inequality, we easily conclude that \eqref{owij} indeed identifies a bounded linear operator in $\mathcal F(\mathcal H)$. In particular, if $f^{(n)}=f_1\otimes\dots\otimes f_n$ with $f_1,\dots,f_n\in B_0(T)$, then
$$ \langle f_1\otimes\dots\otimes f_n,\omega^{\otimes n}\ra=
\langle f_1,\omega\rangle\dotsm\langle f_n,\omega\rangle=x(f_1)\dotsm x(f_n).$$
We will also interpret constants as monomials of order 0.

Let $\mathbf P$  and $\mathbf {CP}$ denote the set of all non-commutative polynomials (finite sums of monomials) with kernels $f^{(n)}\in B_0(T^n)$ and $f^{(n)}\in \D^{(n)}$, respectively.
($\mathbf{CP}$ stands  for ``continuous polynomials.'')
Clearly,
$\mathbf {CP}\subset\mathbf P$ and $\mathbf A\subset\mathbf P$.

\begin{lemma}\label{gtydrt} We have $\mathbf {CP}\Omega=\Ffin$.
\end{lemma}

\noindent {\it Proof}. Clearly, $\mathbf {CP}\Omega\subset\Ffin$. On the other hand, for each $f^{(n)}\in\D^{(n)}$, $$ f^{(n)}=\la f^{(n)},\omega^{\otimes n}\ra\Omega-g^{(n-1)},$$ where $g^{(n-1)}\in\bigoplus_{i=0}^{n-1}\D^{(i)}$. From here, by induction, we conclude that $\Ffin\subset\mathbf {CP}\Omega$.\quad $\square$\vspace{2mm}

We now naturally extend the free expectation $\tau$ to the set $\mathbf{CP}$, and define an inner product
$$(P_1,P_2)_{L^2(\tau)}:=\tau(P_2P_1)=(P_1\Omega,P_2\Omega)_{\mathcal F(\mathcal H)},\quad P_1,P_2\in\mathbf {CP}.
$$

Let $P\in\mathbf {CP}$ and $P\ne0$. Then
$P\ne0$ as an element of $L^2(\tau)$. Indeed, let $P=\sum_{i=0}^n\la f^{(i)},\omega^{\otimes i}\ra$,  where $f^{(n)}\ne0$ (we then call $P$ a polynomial of order $n$). Then the  $\mathcal H^{\otimes n}$-th component of $P\Omega$ is $f^{(n)}$, which implies that $(P,P)_{L^2(\tau)}>0$.
Hence, we can define a real Hilbert space $L^2(\tau)$ as the closure of $\mathbf {CP}$ with respect to the norm generated by the inner product $(\cdot,\cdot)_{L^2(\tau)}$. As we will see below, we can naturally embed $\mathbf P$ (and so also $\mathbf A$) into $L^2(\tau)$. Furthermore, we will also show that every element of $L^2(\tau)$ may be understood as (generally speaking, unbounded) Hermitian operator in $\mathcal F(\mathcal H)$.

Let $\mathbf {CP}^{(n)}$ denote the subset of $\mathbf {CP}$ consisting of all continuous polynomials of order $\le n$.
Let $\mathbf {MP}^{(n)}$ denote the closure of $\mathbf {CP}^{(n)}$ in $L^2(\tau)$.
($\mathbf{MP}$ stands for ``measurable polynomials.'')
Let $\mathbf {OP}^{(n)}:=\mathbf {MP}^{(n)}\ominus \mathbf {MP}^{(n-1)}$, $n\in\mathbb N$, $\mathbf {OP}^{(0)}:=\mathbb R$, where the sign $\ominus$ denotes orthogonal difference in $L^2(\tau)$.
($\mathbf{OP}$ stands for ``orthogonal polynomials.'')
Thus, we get the orthogonal decomposition $L^2(\tau)=\bigoplus_{n=0}^\infty \mathbf {OP}^{(n)}$.

\begin{proposition}\label{yuf}
Consider a linear operator $I:\mathbf {CP}\to\Ffin$ given by $IP=P\Omega$ for $P\in\mathbf {CP}$.
Then,  $I$ extends to a unitary operator $I:L^2(\tau)\to\mathcal F(\mathcal H)$. Furthermore,
$I\mathbf {OP}^{(n)}=\mathcal H^{\otimes n}$.
\end{proposition}

\begin{proof}
The first statement of the proposition directly follows from  Lemma~\ref{gtydrt}. Next,
it follows from the proof of Lemma~\ref{gtydrt} that $I\mathbf {CP}^{(n)}=\bigoplus_{i=0}^n \D^{(i)}$, so that $I\mathbf {MP}^{(n)}=\bigoplus_{i=0}^n{\mathcal H}^{\otimes i}$. From here, the the second statement follows.
\end{proof}

For $f^{(n)}\in\D^{(n)}$, let $P(f^{(n)})$ denote the orthogonal projection of $\langle f^{(n)},\omega^{\otimes n}\rangle$ onto $\mathbf {OP}^{(n)}$, i.e., by the results proved above, $P(f^{(n)})=I^{-1}f^{(n)}$.

\begin{theorem}\label{ftty}
For each $f^{(n)}\in\D^{(n)}$, we have $P(f^{(n)})\in\mathbf {CP}$.
\end{theorem}

Before proving Theorem~\ref{ftty}, we have to introduce some notations. Let $ {NC}(n,\pm1)$ denote the collection
of all \begin{equation}\label{iugytuf}\varkappa=\{(A_1,m_1),\dots(A_k,m_k)\},\quad k\in\mathbb N,\end{equation}
such that $\pi(\varkappa):=\{A_1,\dots,A_k\}$ is an element of ${NC}(n)$, $m_1,\dots,m_k\in\{-1,+1\}$, and if for some $i\in\{1,\dots,k\}$, the set $A_i$ has only one element, then $m_i=1$. For each $j\in\{1,\dots,k\}$, we will interpret $m_j$ as the mark of the element $A_j$ of the non-crossing partition $\pi(\varkappa)$.

Finally, we denote by $G_n$ the subset of
$ {NC}(n,\pm1)$
consisting of all $\varkappa$ as in \eqref{iugytuf} such that there do not exist $i,j\in\{1,\dots,k\}$, $i\ne j$, for which $$ \min A_i<\min A_j\le\max A_j<\max A_i$$ with $m_j=+1$, i.e., an element of a non-crossing partition with mark $+1$ cannot be ``within'' any other element of this partition. (Note that, in \cite{a3,a4}, elements of $G_n$ were called extended partitions, with classes labeled $+1$ called ``classes open on the left''.)

Let $n\in\mathbb N$ and let us fix an arbitrary $\varkappa\in G_n$ as in \eqref{iugytuf}. We then define $W(\varkappa)\omega(t_1)\dotsm\omega(t_n)$ as follows. For each $i\in\{1,\dots,k\}$, let $A_i=\{j_1,j_2,\dots,j_l\}$, $j_1<j_2,\dots<j_l$. If  $m_i=-1$ (and so $l\ge2$), then replace the factors $\omega(t_{j_1}),\omega(t_{j_2}),\dots,\omega(t_{j_l})$ in the product $\omega(t_1)\omega(t_2)\dotsm\omega(t_n)$ by the ``function''
$$\lambda^{l-2}(t_{j_1})\delta(t_{j_1},t_{j_2},\dots,t_{j_l}).$$
If $m_i=+1$, then leave the factor $\omega(t_{j_1})$ without changes, and if $l\ge2$ then additionally replace  the
factors $\omega(t_{j_2}),\omega(t_{j_3}),\dots,\omega(t_{j_l})$ in the product $\omega(t_1)\omega(t_2)\dotsm\omega(t_n)$ by the ``function'
$$\lambda^{l-1}(t_{j_1})\delta(t_{j_1},t_{j_2},\dots,t_{j_l}).$$ Here, analogously to \eqref{hgycdh}, we have set, for $k\ge2$,
$$\int_{T^k}\sigma(t_1)\dotsm\sigma(t_k)\,f^{(k)}(t_1,\dots,t_k)\delta(t_1,t_2,\dots,t_k):=\int_{T}\sigma(dt)\,f^{(k)}(t,t\dots,t).$$

For example, if $n=8$, and
$$\varkappa=\big\{(\{1,2\},+1),(\{3,4,8\},+1),(\{5,6,7\},-1)\big\},$$
then $$W(\varkappa)\omega(t_1)\dotsm\omega(t_8)=
\lambda(t_1)\delta(t_1,t_2)\lambda^2(t_3)\delta(t_3,t_4,t_8)\lambda(t_5)\delta(t_5,t_6,t_7)\omega(t_1)\omega(t_3).$$

Next, we denote by $\operatorname{\it Int}(n)$ the collection of all interval partitions of $\{1,\dots,n\}$, all of whose elements are intervals of consecutive integers. Clearly, $\operatorname{\it Int}(n)\subset {NC}(n)$. We will denote by $\operatorname{\it Int}(n,\pm1)$ the corresponding subset of ${NC}(n,\pm1)$. Note that $\operatorname{\it Int}(n,\pm1)\subset G_n$. 

\begin{proof}[Proof of Theorem~\ref{ftty}]
For any $f\in\D$, denote by $\la f,\omega\ra\cdot$
the operator of left multiplication by $\la f,\omega\ra$ acting on $\mathbf {CP}$. Clearly, under $I$, $\la f,\omega\ra\cdot$ goes over into the operator $\la f,\omega\ra$ acting on $\Ffin$. Now, for any $f_1,\dots,f_n\in\D$, $n\ge2$,
$$ \la f_1,\omega\ra f_2\otimes\dots\otimes f_n=
f_1\otimes\dots\otimes f_n+(\lambda f_1f_2)\otimes f_3\otimes\dots\otimes f_n+(f_1,f_2)_{\mathcal H}f_3\otimes \dots\otimes f_n.$$
Therefore, applying $I^{-1}$ to the above equality, we get
\begin{align}
P(f_1\otimes\dots\otimes f_n)&=\la f_1,\omega\ra P(f_2\otimes\dots\otimes f_n)-P((\lambda f_1f_2)\otimes f_3\otimes\dots\otimes f_n)\notag\\&\quad
-(f_1,f_2)_{\mathcal H}P(f_3\otimes \dots\otimes f_n).\label{tyfsqx}
\end{align}

Let $\D^{(n)}_{\mathrm{alg}}$ denote the subset of $\D^{(n)}$ consisting of finite sums of functions of the form $f_1\otimes\dots\otimes f_n$ with $f_1,\dots,f_n\in\D$. Then, it follows by induction from \eqref{tyfsqx} that, for each $f^{(n)}\in \D^{(n)}_{\mathrm{alg}}$,
\begin{equation}
P(f^{(n)})=\sum_{\varkappa\in \operatorname{\it Int}(n,\pm1)}c(\varkappa)\int_{T^n}\sigma(dt_1)\dotsm
\sigma(t_n)f^{(n)}(t_1,\dots,t_n)W(\varkappa)\omega(t_1)\dotsm\omega(t_n),\label{xffdx}\end{equation}
where $c(\varkappa)\in\mathbb  R$ (compare with \cite[Section~4]{EP} and \cite[Section~3]{a3}).

Now, let us fix an arbitrary $f^{(n)}\in\mathcal D^{(n)}$. Choose a sequence $\{f_k^{(n)}\}_{k=1}^\infty\subset \D^{(n)}_{\mathrm{alg}}$ such that
the set $\bigcup_{k=1}^\infty \operatorname{supp}f_k^{(n)}$ is in $\mathcal B_0(T)$, $f_k^{(n)}$ are uniformly bounded and $f_k^{(n)}\to f^{(n)}$ point-wise as $k\to\infty$. Hence $\la f_k^{(n)},\omega^{\otimes n}\ra\to \la f^{(n)},\omega^{\otimes n}\ra$
in $L^2(\tau)$, which implies that $P(f_k^{(n)})\to P(f^{(n)})$ in $L^2(\tau)$. On the other hand, for each $\varkappa\in G_n$,
\begin{multline*}
\int_{T^n}\sigma(dt_1)\dotsm \sigma(dt_n)\,f_k^{(n)}(t_1,\dots,t_n)W(\varkappa)\omega(t_1)\dotsm\omega(t_n)\\
\to \int_{T^n}\sigma(dt_1)\dotsm \sigma(dt_n)\,f^{(n)}(t_1,\dots,t_n)W(\varkappa)\omega(t_1)\dotsm\omega(t_n)
\end{multline*}
in $L^2(\tau)$ as $n\to\infty$.
This implies that \eqref{xffdx} holds for each $f^{(n)}\in\D^{(n)}$, and therefore $P(f^{(n)})\in\mathbf {CP}$.\end{proof}

For each $n\in\mathbb N$, we define (free) Wick  product of $\omega(t_1),\dots,\omega(t_n)$, denoted by ${:}\omega(t_1)\dotsm\omega(t_n){:}$ as follows: first we formally evaluate the product
$$ \omega(t_1)\dotsm\omega(t_n)=(\di_{t_1}^\dag+\di_{t_1}+\lambda(t_1)\di_{t_1}^\dag\di_{t_1})\dotsm(\di_{t_n}^\dag+\di_{t_n}+\lambda(t_n)\di_{t_n}^\dag\di_{t_n}),$$
and then remove all the terms containing  $\di_{t_{i}}\di_{t_{i+1}}^\dag$ for some $i\in\{1,\dots,n-1\}$. We clearly have the following recursive formula
\begin{align}
&{:}\omega(t_1){:}=\omega(t_1),\notag\\
&{:}\omega(t_1)\dotsm\omega(t_n){:}=\di_{t_1}^\dag{:}\omega(t_2)\dotsm\omega(t_n){:}\notag\\
&\qquad\qquad\qquad\qquad+\lambda(t_1)\di_{t_1}^\dag\di_{t_1}\di_{t_2}\dotsm\di_{t_n}+\di_{t_1}\di_{t_2}\dotsm\di_{t_n},\quad n\ge2.\label{khci}
\end{align}
Furthermore, as easily seen,
\begin{align}
{:}\omega(t_1)\dotsm\omega(t_n){:}&=\di_{t_1}^\dag
\di_{t_2}^\dag\dotsm\di_{t_n}^\dag\notag\\
&\quad +\sum_{i=1}^n(
\di^\dag_{t_1}\dotsm \di_{t_{i-1}}^\dag\di_{t_i}\dotsm \di_{t_n}+\di_{t_1}^\dag\dotsm\di_{t_{i-1}}^\dag\lambda(t_i)\di_{t_i}^\dag\di_{t_i}\di_{t_{i+1}}\dotsm\di_{t_n}).\label{ycfwgy}
\end{align}

\begin{theorem}\label{hjhgs}
For each $f^{(n)}\in\D^{(n)}$, $n\in\mathbb N$,
\begin{equation}
P(f^{(n)})=\int_{T^n}\sigma(dt_1)\dotsm\sigma(dt_n)\,
f^{(n)}(t_1,\dots,t_n)\,{:}\omega(t_1)\dotsm\omega(t_n){:}\,.
\label{hjgscsv}
\end{equation}
\end{theorem}

\begin{proof} Analogously to the proof of Theorem~\ref{ftty}, it suffices to prove formula \eqref {hjgscsv} in the case  $f^{(n)}=f_1\otimes\dots\otimes f_n$ with $f_1,\dots,f_n\in\D$.
Using  \eqref{dhol} and \eqref{khci}, we have:
\begin{align}
&\omega(t_1)\,{:}\omega(t_2)\dotsm\omega(t_n){:}=
\di_{t_1}^\dag\,{:}\omega(t_2)\dotsm\omega(t_n){:}
+(\lambda(t_1)\di_{t_1}^\dag\di_{t_1}+\di_{t_1})(\di_{t_2}^\dag\,{:}\omega(t_3)\dotsm\omega(t_n){:}\notag\\
&\quad+\lambda(t_2)\di_{t_2}^\dag \di_{t_2}\di_{t_3}\dotsm\di_{t_n}+\di_{t_2}\di_{t_3}\dotsm\di_{t_n})\notag\\
&=\di_{t_1}^\dag \,{:}\omega(t_2)\dotsm\omega(t_n){:}+\lambda(t_1)\delta(t_1,t_2)\di_{t_2}^\dag \,{:}
\omega(t_3)\dotsm\omega(t_n){:}\notag\\
&\quad+\lambda(t_1)\delta(t_1,t_2)\lambda(t_2)\di^\dag_{t_2}\di_{t_2}\di_{t_3}\dotsm\di_{t_n}+\lambda(t_1)\di_{t_1}^\dag\di_{t_1}\di_{t_2}\dotsm\di_{t_n}\notag\\
&\quad+\delta(t_1,t_2)\,{:}\omega(t_3)\dotsm\omega(t_n){:}+\lambda(t_1)\delta(t_1,t_2)\di_{t_2}\di_{t_3}\dotsm\di_{t_n}+\di_{t_1}\di_{t_2}\dotsm\di_{t_n}\notag\\
&={:}\omega(t_1)\dotsm\omega(t_n){:}+\lambda(t_1)\delta(t_1,t_2)\,{:}\omega(t_2)\dotsm\omega(t_n){:}+\delta(t_1,t_2)\,{:}\omega(t_3)\dotsm\omega(t_n){:}\,,\label{tyde}
\end{align}
the calculations taking rigorous meaning after smearing out with the $f^{(n)}$ as above.
By virtue of \eqref{tyfsqx}, we see that \eqref{tyde} implies the statement of the theorem.\end{proof}

Taking Theorem~\ref{hjhgs} into account, for each $f^{(n)}\in\D^{(n)}$ we will write $\la f^{(n)},\, {:}\omega^{\otimes n}{:}\ra$ for $P(f^{(n)})$.
More generally, for each $f^{(n)}\in\mathcal H^{\otimes n}$, we will denote by  $\la f^{(n)},\, {:}\omega^{\otimes n}{:}\ra$ the element of $L^2(\tau)$ defined as $I^{-1}f^{(n)}$. Thus, each element $F\in L^2(\tau)$ admits a unique representation
$$ F=\sum_{n=0}^\infty \la f^{(n)},\, {:}\omega^{\otimes n}{:}\ra,$$
where $f=(f^{(n)})\in\mathcal F(\mathcal H)$.

\begin{remark}\label{ghiiuhui}
With each $F\in L^2(\tau)$ one can associate a Hermitian
(i.e., densely defined and symmetric, possibly unbounded) operator in $\mathcal F(\mathcal H)$ with domain $\Ffin$. Indeed, let us fix arbitrary $f^{(n)}\in\D^{(n)}$ and $g^{(m)}\in\D^{(m)}$.
By virtue of \eqref{ycfwgy},
$$ \la f^{(n)},\,{:}\omega^{\otimes n}{:}\ra
g^{(m)}= \mathcal L^{(1)}(f^{(n)})g^{(m)}+\sum_{k=1}^{n\wedge m}\big(\mathcal L_k^{(2)}(f^{(n)})+
\mathcal L_k^{(3)}(f^{(n)})\big)g^{(m)},$$
where
\begin{align}
\mathcal L^{(1)}(f^{(n)})&=\int_{T^n}\sigma(dt_1)\dotsm\sigma(dt_n)\,f^{(n)}(t_1,\dots,t_n)\di_{t_1}^\dag
\dotsm\di_{t_n}^\dag,\notag\\
\mathcal L_k^{(2)}(f^{(n)})&=\int_{T^n}\sigma(dt_1)\dotsm\sigma(dt_n)\,f^{(n)}(t_1,\dots,t_n)\di_{t_1}^\dag\dotsm\di_{t_{n-k}}^\dag\di_{t_{n-k+1}}\dotsm\di_{t_n},\notag\\
\mathcal L_k^{(3)}(f^{(n)})&=\int_{T^n}\sigma(dt_1)\dotsm\sigma(dt_n)\,f^{(n)}(t_1,\dots,t_n)\notag\\
&\qquad\times \di_{t_1}^\dag\dotsm\di_{t_{n-k}}^\dag
\lambda(t_{n-k+1})\di_{t_{n-k+1}}^\dag
\di_{t_{n-k+1}}\di_{t_{n-k+2}}\dotsm\di_{t_n}.\label{jhsdhgttt}
\end{align}
Note that
\begin{gather}
\mathcal L^{(1)}(f^{(n)})g^{(m)}\in\mathcal H^{\otimes(n+m)},\notag \\
\mathcal L_k^{(2)}(f^{(n)})g^{(m)}\in\mathcal H^{\otimes(n+m-2k)},\quad \mathcal L_k^{(3)}(f^{(n)})g^{(m)}\in\mathcal H^{\otimes(n+m-2k+1)}.\label{huddsuy}
\end{gather}
Using \eqref{jhsdhgttt} and the  Cauchy--Schwarz inequality, we conclude that the vectors
in \eqref{huddsuy} are well-defined for each $f^{(n)}\in\mathcal H^{\otimes n}$ (independently of the choice of a version of
$f^{(n)}$), and the $\mathcal F(\mathcal H)$-norm of each such vector is bounded by $C\|f^{(n)}\|_{\mathcal H^{\otimes n}}$,
where the constant $C>0$ only depends on $g^{(m)}$ and is independent of $n$.  Therefore,
for each $F=\sum_{n=0}^\infty\la f^{(n)},\, {:}\omega^{\otimes n}{:}\ra\in L^2(\tau)$,
$$
Fg^{(m)}:=\sum_{n=0}^\infty \mathcal L^{(1)}(f^{(n)})g^{(m)}+\sum_{k=1}^m\bigg(\sum_{n=k}^\infty \big(\mathcal L^{(2)}_k(f^{(n)})+\mathcal L^{(3)}_k(f^{(n)})\big)g^{(m)}\bigg),
$$
which is a vector in $\mathcal F(\mathcal H)$. 
Indeed, by \eqref{huddsuy},
$$
\bigg\|\sum_{n=0}^\infty \mathcal L^{(1)}(f^{(n)})g^{(m)}\bigg\|_{\mathcal F(\mathcal H)}^2=\sum_{n=0}^\infty\|\mathcal L^{(1)}(f^{(n)})g^{(m)}\|_{\mathcal F(\mathcal H)}^2\le C^2 \sum_{n=0}^\infty\| f^{(n)}\|_{\mathcal F(\mathcal H)}^2<\infty,
$$
and analogously we deal with the other sums. 
Extending $F$ by linearity to the whole $\Ffin$, we thus get a Hermitian operator in $\mathcal F(\mathcal H)$
with domain $\Ffin$.
\end{remark}

The following theorem gives a rule of representation of a monomial through a sum of orthogonal polynomials.

\begin{theorem}[Wick rule for a product of  free noises] \label{fyfdytry}
 For each $n\in\mathbb N$, we have:
 \begin{equation}\label{igbyug}
\omega(t_1)\dotsm\omega(t_n)=\sum_{\varkappa\in G_n} {:}W(\varkappa)\omega(t_1)\dotsm\omega(t_n){:}\,,
\end{equation}
the formula making sense after smearing out with a function $f^{(n)}\in\mathcal D^{(n)}$.
\end{theorem}

\begin{proof} We prove \eqref{igbyug} by induction.  Formula  \eqref{igbyug} trivially holds for $n=1$. Assume that it also holds for $n
-1$, $n\ge2$. Then
\begin{align*}
\omega(t_1)\dotsm\omega(t_n)&=\sum_{\varkappa\in G_{n-1}}\omega(t_1)\,{:}W(\varkappa)\omega(t_2)\dotsm\omega(t_n){:}\\
&=\sum_{\varkappa\in G_{n-1}}\sum_{i=1}^3 {:}W(\varkappa^{(i)})\omega(t_1)\dotsm\omega(t_n){:}\, .\end{align*}
Here, $\varkappa^{(i)}$, $i=1,2,3$, are the elements of $ G_n$ that are obtained by  first taking the marked
partition $\varkappa$ of $\{2,3,\dots,n\}$, and then for $i=1$, by adding $\{1\}$ as a singleton element with mark $+1$, for $i=2$, by
adding $1$ to the first (from the left hand side) element of $\varkappa$ which has mark $+1$ (if there is no such an element, then this term is zero),
and for $i=3$,  by adding $1$ to the first element of $\varkappa$ that has mark $+1$ and changing the mark to $-1$ (again this term becomes zero if no element of $\varkappa$ has mark $+1$). From here the statement of the theorem follows.
\end{proof}

\begin{remark}\label{jhgiu}
For each $f^{(n)}\in B_0(T^n)$,
\begin{equation}\label{giuwq} \sum_{\varkappa\in G_n}\int_{T^n}\sigma(dt_1)\dotsm\sigma(dt_n)\,f^{(n)}(t_1,\dots,t_n)\,  {:}W(\varkappa)\omega(t_1)\dotsm\omega(t_n){:}\end{equation}
is clearly an element of $L^2(\tau)$, and it follows from Theorem~\ref{fyfdytry} and Remark~\ref{ghiiuhui} that \eqref{giuwq} is associated with the operator $\la f^{(n)},\omega^{\otimes n}\ra$
(first on $\Ffin$, and then it is extended by continuity to the whole $\mathcal F(\mathcal H)$).
Thus, we get the inclusion of $\mathbf P$ into $L^2(\tau)$.
\end{remark}

The following theorem generalizes Theorem~\ref{fyfdytry}.

\begin{theorem}[Wick rule for a product of normal products of  free noises]\label{gxsresreh}
For any $k_1,\dots,k_l\in\mathbb N$, $l\in\mathbb N$, we have
\begin{align}&
{:}\omega(t_1)\dotsm\omega(t_{k_1}){:}\,{:}\omega(t_{k_1+1})\dotsm \omega(t_{k_1+k_2}){:}\,\dotsm\,
{:}\omega(t_{k_1+k_2+\dots+k_{l-1}+1})\dotsm \omega(t_n){:}\notag\\
&\qquad=\sum {:}W(\varkappa)\omega(t_1)\omega(t_2)\dotsm\omega(t_n){:}\, , \label{drtdserh}
\end{align}
where $n:=k_1+k_2+\dots+k_l$ and the summation in \eqref{drtdserh}
 is over all $\varkappa\in G_n$ such that each element of the induced partition $\pi(\varkappa)$ of
 $\{1,\dots,n\}$ contains maximum one element of  each of the sets
 $\{1,\dots,k_1\}$, $\{k_1+1,\dots,k_1+k_2\}$, \dots, $\{k_1+k_2+\dots+k_{l-1}+1,\dots,n\}$. Formula
 \eqref{drtdserh} makes sense after smearing out with a function $f^{(n)}\in\D^{(n)}$.
\end{theorem}

\begin{proof} Analogously to the proof of Theorem~\ref{fyfdytry}, it suffices to show that, for any
$k_1,k_2\in\mathbb N$,
\begin{align}
&{:}\omega(t_1)\dotsm \omega(t_{k_1}){:}\, {:}\omega(t_{k_1+1})\dotsm\omega(t_{k_1+k_2}){:}\notag\\
&\quad={:}\omega(t_1)\dotsm \omega(t_{k_1}) \omega(t_{k_1+1})\dotsm\omega(t_{k_1+k_2}){:}\notag\\
&\qquad+{:}\omega(t_1)\dotsm \omega(t_{k_1-1}) \delta(t_{k_1},t_{k_1+1}) \omega(t_{k_1+2})
\dotsm\omega(t_{k_1+k_2}){:}\notag\\
&\qquad+{:}\omega(t_1)\dotsm \omega(t_{k_1-1}) \lambda(t_{k_1})\delta(t_{k_1},t_{k_1+1}) \omega(t_{k_1+1})
\dotsm\omega(t_{k_1+k_2}){:}\,.\label{ssou}
\end{align}
To show \eqref{ssou}, represent ${:}\omega(t_1)\dotsm \omega(t_{k_1}){:}$ in the form \eqref{ycfwgy} and represent
\linebreak
$ {:}\omega(t_{k_1+1})\dotsm\omega(t_{k_1+k_2}){:}$ in the form \eqref{khci}, then use
the free commutation relation \eqref{dhol} whenever $\di_{t_{k_1}}\di_{t_{k_1+1}}^\dag$
enters, and finally collect the terms in order to get the right hand side of
\eqref{ssou}. We leave these long, but quite simple calculations to the interested reader.
\end{proof}

\section{ Non-commutative generalized stochastic processes with freely independent values}\label{serts}

Let the space $T$ and the measure $\sigma$ be as in Section~\ref{cfrd}. For each $t\in T$, let $\mu(t,\cdot)$ be a probability measure on $(\mathbb R,\mathcal B (\mathbb R))$ with compact support.  We will assume that, for each $A\in\mathcal B(\mathbb R)$, the mapping $T\ni t\mapsto \mu(t,A)$ is measurable, and for each $\Delta\in\mathcal B_0(T)$ there exists $R=R(\Delta)>0$ such that, for all $t\in\Delta$, the measure $\mu(t,\cdot)$ has support in $[-R,R]$. We denote  $\tilde T:=T\times\mathbb R$, and define a measure
 $\tilde\sigma(dt,ds):=\sigma(dt)\mu(t,ds)$  on $(\tilde T,\mathcal B(\tilde T))$.
Clearly,
\begin{equation}\label{uftufoi}
\tilde\sigma(\Delta\times\mathbb R)<\infty\quad\text{for all }\Delta\in\mathcal B_0(T).\end{equation}
 We denote
$ {\mathcal H}:=L^2(\tilde T,\tilde\sigma)$.
Let  $ \lambda\in C(\tilde T)$ be chosen as $\lambda(t,s):=s$.
Let $L^2(\tau)$ be the  Hilbert space as in Section~\ref{cfrd} which corresponds to $\tilde T$, $\tilde\sigma$, and $\lambda$. By Proposition~\ref{yuf},
we have a unitary operator $I: L^2(\tau)\to \mathcal F({\mathcal H})$.

\begin{remark}
In view of \eqref{uftufoi}, we will call a subset  of $\tilde T$ bounded if it is a subset of a set $\Delta\times \mathbb R$, where $\Delta\in\mathcal B_0(T)$. We then define $B_0(\tilde T)$ and $C_0(\tilde T)$ as the set of all bounded measurable functions on $\tilde T$ with bounded support, and
the set of all bounded continuous functions on $\tilde T$ with bounded support, respectively. All the respective definitions and results of Section~\ref{cfrd} evidently remain true for these spaces.

For each $f:T\to\mathbb R$ and $g:\mathbb R\to\mathbb R$, we denote by $f\otimes g$ the function on $\tilde T$  given by
$(f\otimes g)(t,s):=f(t)g(s)$. If $f\in\D=C_0(T)$ and $g$ is continuous, then $f\otimes g$ is continuous, has bounded support, but is not necessarily bounded. Still we will identify this function with any $f\otimes \bar g\in C_0(\tilde T)$, where $\bar g:\mathbb R\to\mathbb R$ is continuous, bounded, and coincides with $g$ on $[-R,R]$. Here $R=R(\operatorname{supp}f)>0$, i.e., $R$ is chosen so that, for each $t$ from the support of $f$, $\mu(t,\cdot)$ has support in $[-R,R]$.  We will  analogously proceed in the case where $f\in B_0(T)$. \end{remark}

Now, for each $f\in B_0(T)$ we define $X(f)$ as the element of $L^2(\tau)$ given by 
\begin{equation}\label{urytrf}X(f):=x(f\otimes1)=a^+(f\otimes 1)+a^-(f\otimes1)+a^0(f\otimes s).\end{equation}
(Here and below, if  $g(s)=s^l$, $l\in\mathbb N_0$,
we write the function $f\otimes g$ as $f\otimes s^l$.)
Thus,
\begin{align}X(f)&=\int_{\tilde T}\tilde\sigma(dt,ds)f(t)(\di_{(t,s)}^\dag+\di_{(t,s)}+s\di_{(t,s)}^\dag\di_{(t,s)})\notag\\
&=\int_T\sigma(dt) f(t) \int_{\mathbb R}\mu(t,ds)
(\di_{(t,s)}^\dag+\di_{(t,s)}+s\di_{(t,s)}^\dag\di_{(t,s)})\notag\\
&=\int_T\sigma(dt) f(t) \omega(t),\label{tresas}
\end{align}
where
\begin{equation}\label{ertsk}\omega(t):=\int_{\mathbb R}\mu(t,ds)\,
(\di_{(t,s)}^\dag+\di_{(t,s)}+s\di_{(t,s)}^\dag\di_{(t,s)})=\int_{\mathbb R}\mu(t,ds)\varpi(t,s)\end{equation}
with
\begin{equation}\label{jftyfytkkjh} \varpi(t,s)=
\di_{(t,s)}^\dag+\di_{(t,s)}+s\di_{(t,s)}^\dag\di_{(t,s)}.\end{equation}
Also for $\Delta\in\mathcal B_0(T)$, we set $X(\Delta):=X(\chi_\Delta)$.

By Proposition~\ref{dreshg}, for any $f_1,\dots,f_n\in B_0(T)$ such that $f_if_j=0$ $\sigma$-a.e. for all $1\le i< j\le n$, $X(f_1),\dots,X(f_n)$  are freely independent with respect to the state $\tau$. In particular, for any mutually disjoint sets $\Delta_1,\dots,\Delta_n\in\mathcal B_0(T)$, the operators $X(\Delta_1),\dots,X(\Delta_n)$ are freely independent.
Hence, we may interpret $\omega$ as a non-commutative generalized stochastic process with freely independent values (compare with \cite{biane}).

\begin{remark}\label{sreseys}
Let us derive an equivalent representation of the free random field \linebreak $(X(f))_{f\in B_0(T)}$.
For each $t\in T$, denote $c(t):=\mu(t,\{0\})$, and let $\nu(t,\cdot)$ denote the measure on $\mathbb R\setminus\{0\}$ given by $\nu(t,ds):=\frac1{s^2}\,\mu(t,ds)$.
Then, we  define a unitary operator
$$ U:{\mathcal H}\to L^2(T,c(t)\sigma(dt))\oplus
L^2(T\times(\mathbb R\setminus\{0\}),\sigma(dt)\nu(t,ds)) :=\mathcal G$$
by
$${\mathcal H}\ni f\mapsto Uf:=(f(t,0),f(t,s)s)\in \mathcal G. $$
We  naturally extend $U$ to a unitary operator
$U:\mathcal F({\mathcal H})\to\mathcal F(\mathcal G)$. As easily seen, for each $f\in B_0(T)$,
\begin{equation}\label{yddsrts} UX(f)U^{-1}=a^+(f,0)+a^-(f,0)+
a^+(0,f\otimes s)+a^0(0,f\otimes s)+a^-(0,f\otimes s).\end{equation}
In \eqref{yddsrts}, the operator $$B(f):=a^+(f,0)+a^-(f,0)$$ describes the Brownian  part of the process,
while the operator
\begin{equation}\label{esresgy} J(f):=a^+(0,f\otimes s)+a^0(0,f\otimes s)+a^-(0,f\otimes s) \end{equation}
describes the ``jump'' part of the process. Thus,  $\nu(t,\cdot)$ is the L\'evy measure of the process at point $t$, and it describes the value and intensity of ``jumps'' (compare with e.g.\ \cite{parthasarathy} in the bosonic (classical) case, and with \cite{bnt1,bnt2,bnt3} in the free case).
\end{remark}

Analogously to Section~\ref{cfrd}, we define the free cumulants $ C^{(n)}:B_0(T)_{\mathbb C}^n\to\mathbb C$ through
$$ \tau(X(f_1)X(f_2)\dotsm X(f_n))=\sum_{\pi\in {NC}(n)} \prod_{A\in\pi}C(A,f_1,\dots,f_n),\quad f_1,f_2,\dots,f_n\in B_0(T),$$
and then we define the free cumulant transform
$$ C(f):=\sum_{n=1}^\infty C^{(n)}(f),\quad f\in B_0(T)_{\mathbb C}$$
(we have used  obvious notations).
By (the proof of) Proposition~\ref{yuftdfrl} and using the notations introduced in Remark~\ref{sreseys}, we get:

\begin{proposition}\label{hdrtdx}
Let $f\in B_0(T)_{\mathbb C}$ be such that there exists $\varepsilon\in(0,1)$ for which
$$ |f(t)|<\frac{1-\varepsilon}{R}\quad\text{for all }t\in T,$$
where $R=R(\operatorname{supp}f)>0$, i.e., $R$ is such that, for each $t\in\operatorname{supp} f$, the measure $\mu(t,\cdot)$ has support in $[-R,R]$.
Then
\begin{align*}  C (f)&=\int_T\sigma(dt)\int_{\mathbb R}
\mu(t,ds)\,\frac{f^2(t)}{1-sf(t)}\\
&=\int_T \sigma(dt) c(t) f^2(t)+\int_T\sigma(dt)\int_{\mathbb R\setminus\{0\}}
\nu(t,ds)\,\frac{f^2(t)s^2}{1-sf(t)}\\
&=\int_T \sigma(dt) c(t) f^2(t)+\int_T\sigma(dt)\int_{\mathbb R\setminus\{0\}}
\nu(t,ds)\sum_{n=2}^\infty s^nf^n(t).\end{align*}

\end{proposition}

Next, we have:

\begin{proposition}\label{ydtfr}
The vacuum vector $\Omega$ in $\mathcal F({\mathcal H})$ is cyclic for the operator family $(X(f))_{f\in\mathcal D}$. \end{proposition}

\begin{proof}
It can be easily shown by approximation that it suffices to prove that $\Omega$ is cyclic for the operator family $(X(f))_{f\in B_0(T)}$.

 We first state that  the linear span of the set $$\{\chi_\Delta\otimes s^n\mid \Delta\in\mathcal B_0(T),\, n\in\mathbb N_0\}$$ is dense in $L^2(\tilde T,\tilde\sigma)$. Indeed, let $g\in L^2(\tilde T,\tilde\sigma)$ be orthogonal to all elements of this set, i.e.,
\begin{equation}\label{rtsre} \int_\Delta\sigma(dt)\int_{\mathbb R}\mu(t,ds) s^n g(t,s)=0\quad \text{for all $\Delta\in\mathcal B_0(T)$ and $n\in\mathbb N_0$.}\end{equation} Since
$\int_{\mathbb R}\mu(\cdot,ds)s^ng(\cdot,s)\in L^1(\Delta,\sigma)$ for each $\Delta\in \mathcal B_0(T)$, we conclude from \eqref{rtsre} that, for $\sigma$-a.e.\ $t\in T$,
$$ \int_{\mathbb R}\mu(t,ds)s^ng(t,s)=0\quad\text{for all }n\in\mathbb N_0.$$
But, for each $t\in T$, $\mu(t,\cdot)$ is a probability measure on $\mathbb R$ with compact support, and hence the set of all polynomials on $\mathbb R$ is dense in $L^2(\mathbb R,\mu(t,\cdot))$.
Therefore, for $\sigma$-a.e.\ $t
\in T$ and for $\mu(t,\cdot)$-a.e.\ $s\in \mathbb R$, $g(t,s)=0$. Hence $g(t,s)=0$ for $\tilde \sigma$-a.e.\ $(t,s)\in \tilde T$.

Since the measure $\sigma$ is non-atomic, we can analogously prove the following lemma.

\begin{lemma}\label{guygu}
For each $n\in\mathbb N$,
\begin{multline*} {\mathcal H}^{\otimes n}=\operatorname{c.l.s.}(\big\{\chi_{\Delta_1}\otimes s^{l_1})\otimes(\chi_{\Delta_2}\otimes s^{l_2})\otimes\dotsm \otimes (\chi_{\Delta_n}\otimes s^{l_n})\mid
l_1,\dots,l_n\in\mathbb N_0,\\ \text{
for each $j=1,\dots,n-1$: $\Delta_j\cap\Delta_{j+1}=\varnothing$}\big\}.
\end{multline*}
\end{lemma}

Below, we denote by $M$ the set of all multi-indices of the form $(l_1,\dots,l_i)\in\mathbb N_0^i$, $i\in\mathbb N$.

\begin{lemma}\label{huui} For each $n\in\mathbb N$, we define the following subsets of $\mathcal F({\mathcal H})$:
\begin{align}
&\mathcal R^{(n)}:=\operatorname{c.l.s.}\big\{\Omega,\ X(f_1)\dotsm X(f_i)\Omega\mid
f_1,\dots,f_i\in B_0(T),\
i\in\{1,\dots,n\}
\big\},\label{rdsrew}\\
&\mathcal S^{(n)}:= \operatorname{c.l.s.}\big\{\Omega,\ (\chi_{\Delta_1}\otimes s^{l_1})\otimes\dotsm\otimes(\chi_{\Delta_i}\otimes s^{l_i})\mid
 (l_1,\dots,l_i)\in M,\notag\\
 &\qquad l_1+\dots+l_i+i\le n,\
 \text{for each $j=1,\dots,i-1$:}\ \Delta_j\cap\Delta_{j+1}=\varnothing
 \big\}.\notag
\end{align}
Then $\mathcal R^{(n)}=\mathcal S^{(n)}$.
\end{lemma}

\begin{proof}
First, we note by approximation that, for each $n\in\mathbb N$,
\begin{multline}
\mathcal S^{(n)}=\operatorname{c.l.s.}\big\{\Omega,\
f^{(i)}(t_1,\dots,t_i)s_1^{l_1}\dotsm s_i^{l_i}\mid
 f^{(i)}\in B_0(T^i),\\
 (l_1,\dots,l_i)\in M,\ l_1+\dots+l_i+i\le n
\big\}\label{gsdvc}
\end{multline}
(we are using obvious notations for elements of $\mathcal F({\mathcal H})$). From \eqref{urytrf} and \eqref{gsdvc}, the inclusion $\mathcal R^{(n)}\subset\mathcal S^{(n)}$ follows by induction.

Next, let us prove that $\mathcal S^{(n)}\subset\mathcal R^{(n)}$. For $n=1$, this is trivially true. Assume now that this is true for $n\in\{1,\dots,N\}$, and let us show it for $n=N+1$. Thus, we have to show that, for for any
$\Delta_1,\dots,\Delta_i\in\mathcal B_0(T)$ such that $\Delta_j\cap\Delta_{j+1}=\varnothing$ for all
$j=1,\dots,i-1$, and any $(l_1,\dots,l_i)\in M$ such that $l_1+\dots+l_i+i=N+1$,
\begin{equation}\label{ftydyt} (\chi_{\Delta_1}\otimes s^{l_1})\otimes\dotsm \otimes (\chi_{\Delta_i}\otimes s^{l_i})\in\mathcal R^{(N+1)}.\end{equation}

 If $l_1=0$, then
\begin{align}
&(\chi_{\Delta_1}\otimes1)\otimes(\chi_{\Delta_2}\otimes s^{l_2})\otimes \dotsm \otimes(\chi_{\Delta_i}\otimes s^{l_i})\notag\\
&\quad =a^+(\chi_{\Delta_1}\otimes 1)\big( (\chi_{\Delta_2}\otimes s^{l_2})\otimes\dotsm\otimes (\chi_{\Delta_i} \otimes s^{l_i})\big)
\notag\\
&\quad =X(\Delta_1)
\big((\chi_{\Delta_2}\otimes s^{l_2})\otimes\dotsm \otimes(\chi_{\Delta_i}\otimes s^{l_i})\big).
\notag
\end{align}
Hence, in the case  $l_1=0$, \eqref{ftydyt} holds.

Now, for $l_1\ge1$,  we have:
\begin{align}
&(\chi_{\Delta_1}\otimes s^{l_1})\otimes\dotsm\otimes (\chi_{\Delta_i}\otimes s^{l_i})\notag\\
&\quad =a^0(\chi_{\Delta_1}\otimes s)\big((\chi_{\Delta_1}\otimes s^{l_1-1})\otimes (\chi_{\Delta_2}\otimes s^{l_2})\otimes\dotsm\otimes(\chi_{\Delta_i}\otimes s^{l_i})\big)\notag\\
&\quad = X(\Delta_1)\big((\chi_{\Delta_1}\otimes s^{l_1-1})\otimes
(\chi_{\Delta_2}\otimes s^{l_2})\otimes \dotsm\otimes (\chi_{\Delta_i}\otimes s^{l_i})\big)\notag\\
&\qquad\text{} - (\chi_{\Delta_1}\otimes1)\otimes (\chi_{\Delta_1}\otimes s^{l_1-1})\otimes
(\chi_{\Delta_2}\otimes s^{l_2})\otimes \dotsm\otimes (\chi_{\Delta_i}\otimes s^{l_i})\notag\\
&\qquad\text{}-\bigg(\int_{\tilde T}\tilde\sigma(dt,ds)\chi_{\Delta_1}(t)s^{l_1-1}\bigg)(\chi_{\Delta_2}\otimes s^{l_2})\otimes \dotsm\otimes (\chi_{\Delta_i}\otimes s^{l_i}).
\notag
\end{align}
By the results proved above and the induction's assumption,
we therefore conclude that \eqref{ftydyt} holds for $l_1\ge1$.
\end{proof}

From Lemmas \ref {guygu} and \ref{huui} the proposition follows.
\end{proof}

For each $f^{(n)}\in B_0(T^n)$, we define a monomial of $\omega$ by
\begin{align*}
\langle f^{(n)},{\omega}^{\otimes n}\ra:=&\int_{T^n}
\sigma(dt_1)\dotsm\sigma(dt_n)f^{(n)}(t_1,\dots,t_n)
\omega_1(t_1)\dotsm \omega_n(t_n)
\\ =& \int_{\tilde T^n}\tilde\sigma(dt_1,ds_1)\dotsm
\tilde\sigma(dt_n,ds_n)f^{(n)}(t_1,\dots,t_n)\varpi(t_1,s_1)\dotsm \varpi(t_n,s_n)
\end{align*}
(recall \eqref{tresas}--\eqref{jftyfytkkjh}).
We clearly have, for
$f^{(n)}=f_1\otimes\dots\otimes f_n$ with $f_1,\dots,f_n\in B_0(T)$:
\begin{equation}\label{yufytdgg} \langle f_1\otimes\dots\otimes f_n,\omega^{\otimes n}\ra=
\langle f_1,\omega\rangle\dotsm\langle f_n,\omega\rangle=X(f_1)\dotsm X(f_n).\end{equation}
With some abuse of notations, we will  denote by ${\mathbf P}$ and $\mathbf{C}{\mathbf P}$  the set of all polynomials in $\omega$ with kernels $f^{(n)}\in B_0(T^n)$ and $f^{(n)}\in \D^{(n)}$, respectively.
(Note that below we will not use polynomials in the $\varpi$ variable, so keeping the same notations as in  Section~\ref{cfrd} for rather different objects should not lead to a contradiction, and will be justified below.) From Proposition~\ref{ydtfr}, we now conclude:

\begin{proposition}\label{fdestsw5}
The set $\mathbf{C}{\mathbf P}$ is dense in $L^2(\tau)$.
\end{proposition}

Let $\mathbf{C}{\mathbf P}^{(n)}$ denote the subset of $\mathbf{C}{\mathbf P}$ consisting of all continuous polynomials in $\omega$ of order $\le n$.
Let $\mathbf{M}{\mathbf P}^{(n)}$ denote the closure of $\mathbf{C}{\mathbf P}^{(n)}$ in $L^2(\tau)$.
Let $\mathbf{O}{\mathbf P}^{(n)}:=\mathbf{M}{\mathbf P}^{(n)}\ominus \mathbf{M}{\mathbf P}^{(n-1)}$, $n\in\mathbb N$, $\mathbf O{\mathbf P}^{(0)}:=\mathbb R$. Thus, we get:

 \begin{theorem}\label{see3wq}
 We have the following orthogonal decomposition of $L^2(\tau)$:
 $$L^2(\tau)=\bigoplus_{n=0}^\infty \mathbf O{\mathbf P}^{(n)}.$$
 \end{theorem}

Let us recall that, in the case of a classical L\'evy process, Nualart and Schoutens \cite{NS} derived an orthogonal decomposition of any square-integrable functional of the process in multiple stochastic integrals 
with respect to orthogonalized power jump processes (see also  \cite{Ly2} and \cite{a4} for extensions of this result). Our next aim is to derive  a free counterpart of \cite{NS, Ly2}.

Fix any $t\in T$.  Denote by
$(p^{(n)}(t,\cdot))_{n\ge0}$ the system of monic polynomials on $\mathbb R$ which are orthogonal with respect to $\mu(t,\cdot)$. If the support of $\mu(t,\cdot)$ is an infinite set, then by Favard's theorem, the following recursive formula holds:
\begin{align}s p^{(0)}(t,s)&=p^{(1)}(t,s)+b^{(0)}(t),\notag\\
s p^{(n)}(t,s)&=p^{(n+1)}(t,s)+b^{(n)}(t)p^{(n)}(t,s)+a^{(n)}(t)p^{(n-1)}(t,s),\quad n\in\mathbb N,
\label{ydrtyt}\end{align}
where $p^{(0)}(t,s)=1$, $a^{(n)}(t)>0$ for  $n\in\mathbb N$, and $b^{(n)}(t)\in\mathbb R$ for $n\in\mathbb N_0$.
If, however, the support of $\mu(t,\cdot)$ is a finite set consisting of $N$ points ($N\in\mathbb N$), then we have a finite system of monic orthogonal polynomials
$(p^{(n)}(t,\cdot))_{n=0}^{N-1}$ satisfying \eqref{ydrtyt} for $n\le N-2$, and, for $n=N-1$, we have:
$$ sp^{(N-1)}(t,s)=b^{(N-1)}(t,s)p^{(N-1)}(t,s)+a^{(N-1)}(t,s)p^{(N-2)}(t,s).$$ For technical reasons, we  set, in this case,
$$ p^{(n)}(t,s):=0,\  a^{(n)}(t):=0,\quad n\ge N,$$
($b^{(n)}(t)$, $n\ge N$ being arbitrary), so that recursive relation \eqref{ydrtyt} now always holds.

For each $n\in \mathbb N_0$, we denote
\begin{equation}\label{tdtdtdftt} g^{(l)}(t):=\int_{\mathbb R}\mu(t,ds)|p^{(l)}(t,s)|^2,\quad t\in T,\end{equation}
and then we define a measure on $(T,\mathcal B(T))$ by
\begin{equation}\label{sresydrtdstrs} \sigma^{(l)}(dt):=g^{(l)}(t)\sigma(dt).\end{equation}
Note that $\sigma^{(0)}=\sigma$.
 For each $(l_1,\dots,l_i)\in M$, we define
\begin{equation}\label{rtsryhgy} \mathbb H_{(l_1,\dots,l_i)}:=L^2(T^i,\sigma^{(l_1)}\otimes\dots\otimes \sigma^{(l_i)}).\end{equation}
Then, clearly, the following mapping is an isometry
\begin{multline}\label{yurftydeff}\mathbb H_{(l_1,\dots,l_i)}\ni f^{(i)}\mapsto K_{(l_1,\dots,l_i)}f^{(i)}=
(K_{(l_1,\dots,l_i)}f^{(i)})(t_1,s_1,\dots,t_i,s_i)\\
:=
f^{(i)}(t_1,\dots,t_i)p^{(l_1)}(t_1,s_1)\dotsm p^{(l_i)}(t_i,s_i)\in{\mathcal H}^{\otimes i}.
\end{multline}
We denote by $\mathcal H_{(l_1,\dots,l_i)}$ the range of the isometry $K_{(l_1,\dots,l_i)}$.

\begin{lemma}\label{igtyufcf} We have
\begin{equation}\label{cdrtdr} \mathcal F( H)=\mathbb R\oplus \bigoplus_{(l_1,\dots,l_i)\in M}\mathcal H_{(l_1,\dots,l_i)}.\end{equation}
Furthermore,
for each $(l_1,\dots,l_i)\in M$, we have:
\begin{multline}
\mathcal H_{(l_1,\dots,l_i)}=\operatorname{c.l.s.}
\big\{
(\chi_{\Delta_1}\times p^{(l_1)})\otimes\dots\otimes(\chi_{\Delta_i}\times p^{(l_i)})\mid\\
\text{$\Delta_1,\dots,\Delta_i\in\mathcal B_0(T)$},\
\text{for all $j=1,\dots,i-1$:}\  \Delta_j\cap\Delta_{j+1}=\varnothing
\big\}.\label{cfdydwqwqqw}
\end{multline}
Here, $(\chi_{\Delta}\times p^{(l)})(t,s):=\chi_{\Delta}(t)p^{(l)}(t,s)$.
\end{lemma}

\begin{proof} Fix any $f^{(i)}\in B_0(T^i)$. Let $\Delta\in\mathcal B_0(T)$ be such that the support of $f^{(i)}$ is a subset of $\Delta^i$. Choose $R=R(\Delta)>0$ such that, for each $t\in\Delta$, $\mu(t,\cdot)$ has support in $[-R,R]$. Recall the recursive formula \eqref{ydrtyt}. We have, for each $t\in\Delta$, $|b^{(n)}(t)|\le R$ and  $a^{(n)}(t)\le R^2$,
which easily follows from the theory of Jacobi matrices (see e.g.\ \cite{ber}). Therefore, by \eqref{ydrtyt}, each $p^{(n)}(t,s)$ is bounded as a function of $(t,s)\in\Delta\times [-R,R]$. Therefore, for each $f^{(i)}\in B_0(T^i)$,
$$ f^{(i)}(t_1,\dots,t_i)p^{(l_1)}(t_1,s_1)\dotsm p^{(l_i)}(t_i,s_i)\in \mathcal H_{(l_1,\dots,l_i)}.$$
From here equality \eqref{cfdydwqwqqw} easily follows (recall that the measure $\sigma$ is non-atomic, which allows us to choose only those sets $\Delta_1,\dots,\Delta_i$ in \eqref{cfdydwqwqqw} for which $\Delta_j\cap\Delta_{j+1}=\varnothing$ for $j=1,\dots,i-1$). Formula \eqref{cdrtdr} can now be  proven analogously to the proof of Lemma~\ref{guygu}.
\end{proof}

Recall that by Proposition~\ref{yuf}, we have a unitary operator $I:L^2(\tau)\to\mathcal F({\mathcal H})$. For each $(l_1,\dots,l_i)\in M$, denote $\mathbf H_{(l_1,\dots,l_i)}.:=I^{-1}\mathcal H_{(l_1,\dots,l_i)}$.
For any $\Delta\in\mathcal B_0(T)$ and $l\in\mathbb N_0$, denote
$$ X^{(l)}(\Delta):= \int_{\tilde T}\tilde \sigma(dt,ds) \chi_\Delta (t) p^{(l)}(t,s)\varpi(t,s).$$
For arbitrary $(l_1,\dots,l_i)\in M$ and $\Delta_1,\dots,\Delta_i\in\mathcal B_0(T)$ such that $\Delta_j\cap\Delta_{j+1}=\varnothing$ for $j=1,\dots,i-1$,  we clearly have:
$$X^{(l_1)}(\Delta_1)\dotsm X^{(l_i)}(\Delta_i)\Omega=
(\chi_{\Delta_1}\times p^{(l_1)})\otimes\dots\otimes(\chi_{\Delta_i}\times p^{(l_i)}).
$$

 Therefore, by \eqref{cfdydwqwqqw},
\begin{multline*}
\mathbf H_{(l_1,\dots,l_i)}=\operatorname{c.l.s.}
\big\{ X^{(l_1)}(\Delta_1)\dotsm X^{(l_i)}(\Delta_i)
\mid
\text{$\Delta_1,\dots,\Delta_i\in\mathcal B_0(T)$},
\\
\text{for all $j=1,\dots,i-1$:}\  \Delta_j\cap\Delta_{j+1}=\varnothing
\big\}.
\end{multline*}

For each $f^{(l_1,\dots,l_i)}\in \mathbb H_{(l_1,\dots,l_i)}$ (recall \eqref{rtsryhgy}), we can easily
define a non-commutative multiple stochastic integral
\begin{equation}\label{tsreseras}
\int_{T^i} f^{(l_1,\dots,l_i)}(t_1,\dots,t_i)X^{(l_1)}(dt_1)\dotsm X^{(l_i)}(dt_i)
\end{equation}
as an element of $\mathbf H_{(l_1,\dots,l_i)}$. Indeed, for each $f^{(l_1,\dots,l_i)}$ of the form
$$ f^{(l_1,\dots,l_i)}(t_1,\dots,t_i)=\chi_{\Delta_1}(t_1)\dotsm \chi_{\Delta_i}(t_i)$$
with $\Delta_1,\dots,\Delta_i\in\mathcal B_0(T)$ such that $\Delta_j\cap\Delta_{j+1}=\varnothing$, $j=1,\dots,i-1$, we define
\eqref{tsreseras} as $X^{(l_1)}(\Delta_1)\dotsm X^{(l_i)}(\Delta_i)$.
We then extend this definition by linearity to the linear span of such functions,
and finally we extend it by continuity to obtain a unitary operator
$$ \mathbb H_{(l_1,\dots,l_i)}\ni f^{(l_1,\dots,l_i)}\mapsto \int_{T^i} f^{(l_1,\dots,l_i)}(t_1,\dots,t_i)X^{(l_1)}(dt_1)\dotsm X^{(l_i)}(dt_i)\in  \mathbf H_{(l_1,\dots,l_i)}.$$

Taking \eqref{cdrtdr} into account, we thus derive

\begin{theorem}\label{fytytif}
Denote
$$ \mathbb F:= \mathbb R\oplus \bigoplus_{(l_1,\dots,l_i)\in M}\mathbb H_{(l_1,\dots,l_i)}.$$
Then,  the following unitary operator gives an orthogonal expansion of $L^2(\tau)$ in non-commutative multiple stochastic integrals:
\begin{multline}\label{vfytfytd}
\mathbb F\ni F=(c, (f^{(l_1,\dots,l_i)})_{(l_1,\dots,l_i)\in M})\\
\mapsto JF:=
c\mathbf1+\sum_{(l_1,\dots,l_i)\in M}\int_{T^i} f^{(l_1,\dots,l_i)}(t_1,\dots,t_i)X^{(l_1)}(dt_1)\dotsm X^{(l_i)}(dt_i)\in L^2(\tau).
\end{multline}
In terms of this orthogonal expansion, we have:
\begin{equation}\label{hguyg} L^2(\tau)=\mathbb R\oplus \bigoplus_{(l_1,\dots,l_i)\in M}\mathbf H_{(l_1,\dots,l_i)}.\end{equation}
(Note that, in \eqref{hguyg}, $\mathbb R$ denotes the space of all operators $c\mathbf 1$, where $c\in\mathbb R$.)

\end{theorem}

\begin{remark} For each $l\in \mathbb N_0$ and $\Delta\in\mathcal B_0(T)$, define $Y^{(l)}(\Delta)\in L^2(\tau)$ by
\begin{align*}
Y^{(l)}(\Delta):&=\int_{\tilde T}\tilde\sigma (dt,ds)\chi_\Delta(t)s^l\varpi(t,s)\\
&=a^+(\chi_\Delta\otimes s^l)+a^0(\chi_\Delta\otimes s^{l+1})+a^-(\chi_\Delta\otimes s^l)
\end{align*}
(recall \eqref{jftyfytkkjh}).
Clearly, $Y^{(0)}(\Delta)=X(\Delta)$. Recall now the unitary operator $U:\mathcal F({\mathcal H})\to\mathcal F(\mathcal G)$ from Remark~\ref{sreseys}. Then, for each $l\in\mathbb N$, we have:
$$ U Y^{(l)}(\Delta)U^{-1}=a^{+}(0,\chi_\Delta\otimes s^{l+1}) +a^{0}(0,\chi_\Delta\otimes s^{l+1})+a^{-}(0,\chi_\Delta\otimes s^{l+1})$$
(compare with \eqref{yddsrts} and \eqref{esresgy}). Hence, by analogy with the classical case (see \cite{NS}), $Y^{(l)}(\cdot)$, $l\in\mathbb N_0$, may be treated
as ``power jump processes'' (recall that $s$ describes the value of ``jumps'').  For any $l_1,l_2\in\mathbb N_0$, $l_1<l_2$, and any $\Delta_1,\Delta_2\in\mathcal B_0(T)$,
\begin{align*}
\tau(Y^{(l_1)}(\Delta_1)X^{(l_2)}(\Delta_2))&=(X^{(l_2)}(\Delta_2)\Omega,Y^{(l_1)}(\Delta_1)\Omega)_{\mathcal F({\mathcal H})}\\
&=\int_{\Delta_1\cap\Delta_2}\sigma(dt)\int_{\mathbb R}\mu(t,ds)p^{(l_2)}(t,s)s^{l_1}=0.
\end{align*}
Therefore, $X^{(l)}(\cdot)$, $l\in\mathbb N_0$, may be thought of as the orthogonalized power jump processes $Y^{(l)}(\cdot)$, $l\in\mathbb N_0$.
\end{remark}

The following theorem describes a connection between Theorems~\ref{see3wq}  and \ref{fytytif}.

\begin{theorem}\label{drse}
For each $n\in\mathbb N$,
$$ \mathbf O{\mathbf P}^{(n)}=\bigoplus_{(l_1,\dots,l_i)\in M,\  l_1+\dots+l_i+i=n}\mathbf H_{(l_1,\dots,l_i)}.$$
\end{theorem}

\begin{proof} We have to show that, for each $n\in\mathbb N$,
\begin{equation*}
\mathbf M{\mathbf P}^{(n)}=\mathbb R\oplus\bigoplus_{(l_1,\dots,l_i)\in M,\  l_1+\dots+l_i+i\le n}\mathbf H_{(l_1,\dots,l_i)},\end{equation*}
or, equivalently, $I\mathbf M{\mathbf P}^{(n)}=\mathcal Z^{(n)}$, where
\begin{equation}\label{yfctydf}
\mathcal Z^{(n)}:=\mathbb R\oplus\bigoplus_{(l_1,\dots,l_i)\in M,\  l_1+\dots+l_i+i\le n}\mathcal H_{(l_1,\dots,l_i)}.\end{equation}
As easily seen, $I\mathbf M{\mathbf P}^{(n)}=\mathcal R^{(n)}$ (see \eqref{rdsrew}). Hence, by Lemma~\ref{huui}
and \eqref{gsdvc}
\begin{multline}
I\mathbf M{\mathbf P}^{(n)}=\operatorname{c.l.s.}\big\{\Omega,\
f^{(i)}(t_1,\dots,t_i)s_1^{l_1}\dotsm s_i^{l_i}\mid
 f^{(i)}\in B_0(T^i),\\
 (l_1,\dots,l_i)\in M,\ l_1+\dots+l_i+i\le n
\big\}.\label{srehgsdvc}
\end{multline}
Furthermore, \eqref{cfdydwqwqqw} implies that
\begin{multline}
\mathcal Z^{(n)}=\operatorname{c.l.s.}\big\{\Omega,\
f^{(i)}(t_1,\dots,t_i) p^{(l_1)}(t_1,s_1)\dotsm p^{(l_i)}(t_i,s_i)
\mid
 f^{(i)}\in B_0(T^i),\\
 (l_1,\dots,l_i)\in M,\ l_1+\dots+l_i+i\le n
\big\}.\label{yretdstws}
\end{multline}
It follows from the proof of Lemma~\ref{igtyufcf} that each $p^{(l)}$ has a representation
$$ p^{(l)}(t,s)=\sum_{j=0}^l \alpha^{(l,\,j)}(t)s^j, $$
where $\alpha^{(l,\,j)}$'s are measurable functions on $T$ which are bounded on each $\Delta\in\mathcal B_0(T)$. By \eqref{srehgsdvc} and \eqref{yretdstws}, we therefore get the inclusion $\mathcal Z^{(n)}\subset I\mathbf M{\mathbf P}^{(n)}$.

Next, for each $t\in T$, denote by $P^{(i)}(t,\cdot)$, $i\in\mathbb N_0$, the system of {\it normalized} orthogonal polynomials in $L^2(\mathbb R,\mu(t,\cdot))$. We then  have an expansion
\begin{equation}
\label{d6yed65} s^l=\sum_{j=0}^l \beta^{(l,\,j)}(t)P^{(j)}(t,s), \end{equation}
where the  functions $\beta^{(l,\,j)}$ are measurable.
For each $\Delta\in\mathcal B_0(T)$ and each $t\in\Delta$, we have:
$$ \sum_{j=0}^l \big(\beta^{(l,j)}(t)\big)^2=
\|s^l\|^2_{L^2(\mathbb R,\mu(t,ds))}=
\|s^l\|^2_{L^2([-R,R],\mu(t,ds))}\le R^{2l},$$
where $R=R(\Delta)$.
Thus, the functions $\beta^{(l,\,i)}(\cdot)$ are locally bounded on $T$.
Define
\begin{multline*}
\mathcal Y^{(n)}:=\operatorname{c.l.s.}\big\{\Omega,\
f^{(i)}(t_1,\dots,t_i) P^{(l_1)}(t_1,s_1)\dotsm P^{(l_i)}(t_i,s_i)
\mid
 f^{(i)}\in B_0(T^i),\\
 (l_1,\dots,l_i)\in M,\ l_1+\dots+l_i+i\le n
\big\}.
\end{multline*}
Then, by \eqref{srehgsdvc} and \eqref{d6yed65}, $ I\mathbf M{\mathbf P}^{(n)}\subset \mathcal Y^{(n)}$.
Set $u^{(l)}(t):=\|p^{(l)}(t,\cdot)\|^{-1}_{L^2(\mathbb R,\,\mu(t,\cdot))}$.
(In the case where $p^{(l)}(t,\cdot)=0$, set $u^{(l)}(t):=0$.)
 To show that $ I\mathbf M{\mathbf P}^{(n)}\subset \mathcal Z^{(n)}$, it only remains to show that, for each $f^{(i)}\in B_0(T^i)$ and each $(l_1,\dots,l_i)\in M$, $l_1+\dots+l_i+i\le n$, the function
$$ f^{(i)}(t_1,\dots,t_i)u^{(l_1)}(t_1)\dotsm u^{(l_i)}(t_i)p^{(l_1)}(t_1,s_1)\dotsm p^{(l_i)}(t_i,s_i) $$
belongs to $\mathcal Z^{(n)}$. But this easily follows through approximation of $$f^{(i)}(t_1,\dots,t_i)u^{(l_1)}(t_1)\dotsm u^{(l_i)}(t_i)$$ by functions from $B_0(T^i)$.
\end{proof}


Recall that we have constructed the following chain of unitary isomorphisms:
$$ \mathbb F \overset{J}{\longrightarrow} L^2(\tau) \overset{I}{\longrightarrow} \mathcal F({\mathcal H}) $$
(see, in particular, Theorem~\ref{fytytif}). Thus, $$K:=IJ: \mathbb F\to\mathcal F({\mathcal H})$$ is a unitary operator. Note that 
the restriction of $K$ to each space $\mathbb H_{(l_1,\dots,l_i)}$ is $K_{(l_1,\dots,l_i)}$,  see \eqref{yurftydeff}.
We will preserve the notation $\Omega$ for the vector in $\mathbb F$ defined as $K^{-1}\Omega$.

For each $n\in\mathbb N_0$, we denote
$$ \mathbb F^{(n)}:=J^{-1}\mathbf O{\mathbf P}^{(n)},$$
so that $\mathbb F=\bigoplus_{n=0}^\infty \mathbb F^{(n)}$, and by Theorem~\ref{drse}, for each $n\in\mathbb N$,
$$ \mathbb F^{(n)}=\bigoplus_{(l_1,\dots,l_i)\in M,\  l_1+\dots+l_i+i=n}\mathbb H_{(l_1,\dots,l_i)}.$$
For each $f\in B_0(T)$, we will preserve the notation  $X(f)$ for the image of this operator under $K^{-1}$, i.e., for the equivalent realization of $X(f)$ in $\mathbb F$.

\begin{corollary}\label{gyufdyu}
For each $f\in B_0(T)$, we have $X(f)=X^+(f)+X^0(f)+X^-(f)$, where $X^+(f):\mathbb F^{(n)}\to\mathbb F^{(n+1)}$, $X^0(f): \mathbb F^{(n)}\to \mathbb F^{(n)}$, and $X^-(f): \mathbb F^{(n)}\to \mathbb F^{(n-1)}$. Furthermore, $ X^\pm(f)=X_1^\pm(f)+X_2^\pm(f)$, and for each $(l_1,\dots,l_i)\in M$
 \begin{align}&\mathbb H_{(l_1,\dots,l_i)}\ni g\mapsto (X_1^+(f)g)(t_1,\dots,t_{i+1}):=f(t_1)g(t_2,\dots,t_{i+1})\in \mathbb H_{(0,l_1,\dots,l_i)},\label{111}\\
 &\mathbb H_{(l_1,\dots,l_i)}\ni g\mapsto (X_2^+(f)g)(t_1,\dots,t_{i}):=f(t_1)g(t_1,\dots,t_{i})\in \mathbb H_{(l_1+1,l_2,\dots,l_i)},\label{222}\\
 &\mathbb H_{(l_1,\dots,l_i)}\ni g\mapsto (X_1^-(f)g)(t_1,\dots,t_{i-1}) :=\delta_{l_1,\,0}\int_T \sigma(dt)f(t)g(t,t_1,\dots,t_{i-1})\in \mathbb H_{(l_2,\dots,l_i)},\notag\\
  &
  \mathbb H_{(l_1,\dots,l_i)}\ni g\mapsto (X_2^-(f)g)(t_1,\dots,t_{i}):=(1-\delta_{l_1,\,0})a^{(l_1)}(t_1)f(t_1)g(t_1,\dots,t_{i})\in \mathbb H_{(l_1-1,\dots,l_i)},\notag\\
 &\mathbb H_{(l_1,\dots,l_i)}\ni g\mapsto (X^0(f)g)(t_1,\dots,t_{i}):=b^{(l_1)}(t_1)f(t_1)g(t_1,\dots,t_{i})\in \mathbb H_{(l_1,\dots,l_i)}, \notag
 \end{align}
 and $X^0(f)\Omega=X^-(f)\Omega=0$, $X^+(f)\Omega=f\in\mathbb H_{(0)}$. Here, $\delta_{l_1,\,0}$ is equal to $1$ if $l_1=0$, and equal to $0$, otherwise.
\end{corollary}

\begin{proof} We fix any $f\in B_0(T)$ and $g^{(i)}\in\mathbb H_{(l_1,\dots,l_i)}$.
Then, by \eqref{ydrtyt} and \eqref{yurftydeff}, we have:
\begin{align*}
& \big(a^+(f\otimes 1)+a^-(f\otimes 1)+a^0(f\otimes s)\big)Kg^{(i)}\\
&\quad=\big(a^+(f\otimes 1)+a^-(f\otimes 1)+a^0(f\otimes s)\big) g^{(i)}(t_1,\dots,t_i)p^{(l_1)}(t_1,s_1)\dotsm p^{(l_i)}(t_i,s_i)\\
&\quad=f(t_1)g^{(i)}(t_2,\dots,t_{i+1})p^{(0)}(t_1,s_1)p^{(l_1)}(t_2,s_2)\dotsm p^{(l_i)}(t_{i+1},s_{i+1})\\
&\qquad+\delta_{l_1,0}\int_T\sigma(dt)f(t)g^{(i)}(t,t_1,\dots,t_{i-1})p^{(l_2)}(t_1,s_1)\dotsm p^{(l_i)}(t_{i-1},s_{i-1})\\
&\qquad+f(t_1)g^{(i)}(t_1,\dots,t_i)(p^{(l_1+1)}(t_1,s_1)+b^{(l_1)}(t_1)p^{(l_1)}(t_1,s_1)+a^{(l_1)}(t_1)p^{(l_1-1)}(t_1,s_1))\\
&\qquad\quad\times p^{(l_2)}(t_2,s_2)\dotsm p^{(l_i)}(t_i,s_i).
\end{align*}
Applying the operator $K^{-1}$ to the above element of $\mathcal F(\mathcal H)$, we easily conclude the statement. 
\end{proof}

For $f^{(n)}\in\D^{(n)}$, let $ P(f^{(n)})$ denote the orthogonal projection of $\langle f^{(n)},\omega^{\otimes n}\rangle$ onto $\mathbf {O}{\mathbf P}^{(n)}$.

\begin{remark}\label{sra}
By Proposition~\ref{fdestsw5}, the set $\mathbf{C}{\mathbf P}$ is dense in $L^2(\tau)$. From here it follows  that the linear span of the set $\{P(f^{(n)})\mid f^{(n)}\in\D^{(n)},\ n\in\mathbb N_0\}$ is also dense in $L^2(\tau)$. In fact, for each $n\in\mathbb N$, the set  $\{P(f^{(n)})\mid f^{(n)}\in\D^{(n)}\}$ is dense in $\mathbf {O}{\mathbf P}^{(n)}$. Indeed, by definition, the set $\mathbf C{\mathbf P}^{(n)}$ is dense in $\mathbf M{\mathbf P}^{(n)}$. Therefore, the set of all projections of $P\in \mathbf C{\mathbf P}^{(n)}$ onto $\mathbf {O}{\mathbf P}^{(n)}$ is dense in $\mathbf {O}{\mathbf P}^{(n)}$. But the projection of each $P\in \mathbf C{\mathbf P}^{(n-1)}$
onto $\mathbf {O}\mathbf {P}^{(n)}$ equals zero, from where the statement follows. \end{remark}

\begin{corollary}\label{tsreaw} Let $n\in\mathbb N$ and let $(l_1,\dots,l_i)\in M$, $l_1+\dots+l_i+i=n$. For each $f^{(n)}\in\D^{(n)}$, the $\mathbb H_{(l_1,\dots,l_i)}$-coordinate of the vector $J^{-1} P(f^{(n)})$ in $\mathbb F^{(n)}$ is given by
$$ f(\underbrace{t_1,\dots,t_1}_{\text{$(l_1+1)$ times}}, \underbrace{t_2,\dots,t_2}_{\text{$(l_2+1)$ times}},\dots,
\underbrace{t_i,\dots,t_i}_{\text{$(l_i+1)$ times}}). $$
\end{corollary}

\begin{proof}
By approximation, it suffices to check the statement in the case where $f^{(n)}=f_1\otimes\dots\otimes f_n$, $f_1,\dots,f_n\in\D$. Then, by \eqref{yufytdgg},
$$ J^{-1}\langle f^{(n)},\omega^{\otimes n}\rangle= X(f_1)\dotsm X(f_n)\Omega.$$
Hence, $J^{-1} P(f^{(n)})$ is equal to the projection of  $X(f_1)\dotsm X(f_n)\Omega$ onto $\mathbb F^{(n)}$. Therefore, by Corollary~\ref{gyufdyu},
$$ J^{-1} P(f^{(n)})= X^+(f_1)\dotsm X^+(f_n)\Omega.$$
 The statement now follows from \eqref{111} and \eqref{222}.\end{proof}

In view of Remark~\ref{sra} and Corollary~\ref{tsreaw}, we will now give an equivalent interpretation of the  $\mathbb F^{(n)}$ spaces. So, we fix any  $n\in \mathbb N$.  For each $(l_1,\dots,l_i)\in M$, $l_1+\dots+l_i+i=n$, we define
\begin{multline*}
T^{(l_1,\dots,l_i)}:=\big\{
(t_1,\dots,t_n)\in T^n\mid t_1=t_2=\dots=t_{l_1+1},\\ t_{l_1+2}=t_{l_1+3}=\dots=t_{l_1+l_2+2},\dots,
t_{l_1+l_2+\dots+l_{i-1}+i}=t_{l_1+l_2+\dots+l_{i-1}+i+1}=\dots=t_n,\\
t_{l_1+1}\ne t_{l_1+l_2+2},\ t_{l_1+l_2+2}\ne t_{l_1+l_2+l_3+3}\dots, t_{l_1+\dots+l_{i-1}+i-1}\ne t_n
\big\}.
\end{multline*}
The $T^{(l_1,\dots,l_i)}$ sets with $(l_1,\dots,l_i)\in M$, $l_1+\dots+l_i+i=n$, form a set partition of $T^n$.

We define $\mathcal B(T^{(l_1,\dots,l_i)})$ as the trace $\sigma$-algebra of $\mathcal B(T^n)$ on $T^{(l_1,\dots,l_i)}$. Now, consider the measurable mapping
\begin{equation}\label{tydryde}
T^{(l_1,\dots,l_i)}\ni (t_1,\dots,t_n)\mapsto (t_{l_1+1}, t_{l_1+l_2+2},t_{l_1+l_2+l_3+3},\dots, t_n)\in T^i.\end{equation}
Since $\sigma$ is a non-atomic measure, the image of $T^{(l_1,\dots,l_i)}$ under the mapping \eqref{tydryde} is of full $\sigma^{(l_1)}\otimes\dots\otimes \sigma^{(l_i)}$ measure.
We denote by $\gamma^{(l_1,\dots,l_i)}$ the pre-image of the measure $\sigma^{(l_1)}\otimes\dots\otimes \sigma^{(l_i)}$ under the mapping \eqref{tydryde}. We then extend $\gamma^{(l_1,\dots,l_i)}$ by zero to the whole space $T^n$.  Note that, for different $(l_1,\dots,l_i)$ and $(l_1',\dots,l'_{j})$ from $M$ for which $l_1+\dots+l_i+i=l_1'+\dots+l_j'+j=n$, the measures
$\gamma^{(l_1,\dots,l_i)}$ and $\gamma^{(l'_1,\dots,l'_j)}$ are concentrated on disjoint sets in $T^n$. We then define a measure on $(T^n,\mathcal B(T^n))$ as follows:
\begin{equation}\label{11111} \gamma_n:=\sum_{(l_1,\dots,l_i)\in M,\ l_1+\dots+l_i+i=n}\gamma^{(l_1,\dots,l_i)}.\end{equation}

Recall that, by Remark~\ref{sra}, the set
$\{J^{-1}\ P(f^{n})\mid f^{(n)}\in\D^{(n)}\}$ is dense in $\mathbb F^{(n)}$, while the set
$\D^{(n)}$ is clearly dense in $L^2(T^n,\gamma_n)$.
Therefore, by Corollary~\ref{tsreaw} the mapping
$$ L^2(T^n,\gamma_n)\supset \D^{(n)}\ni f^{(n)}\mapsto J^{-1} P(f^{(n)})\in\mathbb F^{(n)}$$
extends to a unitary operator. In terms of this unitary isomorphism, we will, in what follows,
identify $\mathbb F^{(n)}$ with $L^2(T^n,\gamma_n)$, so that the space $\mathbb F$ becomes
 $$
 \mathbb F=\mathbb R\oplus\bigoplus_{n=1}^\infty L^2(T^n,\gamma_n).$$ 
By analogy with \cite{Ly1,Ly2}, we call $\mathbb F$ a free extended Fock space. 
Since, for each $n\in\mathbb N$, $\D^{(n)}\subset L^2(T^n,\gamma_n)$, we have an evident inclusion of $\Ffin$ into $\mathbb F$.  Corollaries~\ref{gyufdyu} and \ref{tsreaw} can now be reformulated as the following theorem, which is the main result of this section.

\begin{theorem}\label{yufytrkjh} The following mapping
\begin{equation}\label{dyrrrs} \mathbb F\supset \Ffin\ni(f^{(0)},f^{(1)},f^{(2)},\dots) \overset{J}{\longrightarrow} \sum_{n=0}^\infty P(f^{(n)})\in L^2(\tau)\end{equation}
(the sum being, in fact, finite) extends to the unitary operator $J:\mathbb F\to L^2(\tau)$. 
In particular,
for any $f^{(n)},g^{(n)}\in\mathcal D^{(n)}$, $n\in\mathbb N$,
\begin{align}
(P(f^{(n)}),P(g^{(n)}))_{L^2(\tau)}&=(f^{(n)},g^{(n)})_{L^2(T^n,\gamma_n)}\notag\\
&=\sum_{(l_1,\dots,l_i)\in M,\ l_1+\dots+l_i+i=n}
\int_{T^i}(f^{(n)}g^{(n)})(\underbrace{t_1,\dots,t_1}_{\text{$l_1+1$ times}},\dots, \underbrace{t_i,\dots,t_i}_{\text{$l_i+1$ times}})\notag\\
&\qquad\qquad\times g^{(l_1)}(t_1)\dotsm g^{(l_i)}(t_i)\,\sigma(dt_1)\dotsm
\sigma(dt_i),\label{iutfrsdr}\end{align}
where the functions $g^{(l)}$ are given by \eqref{tdtdtdftt}.

For each
$f\in \D$,
$X(f)=X^+(f)+X^0(f)+X^-(f)$, where $X^+(f):\mathbb F^{(n)}\to\mathbb F^{(n+1)}$, $X^0(f): \mathbb F^{(n)}\to \mathbb F^{(n)}$, and $X^-(f): \mathbb F^{(n)}\to \mathbb F^{(n-1)}$. Furthermore, for each $n\in\mathbb N$ and each $g^{(n)}\in \D^{(n)}$,
\begin{equation}\label{ydse} (X^+(f)g^{(n)})(t_1,\dots,t_{n+1})=f(t_1)g(t_2,\dots,t_{n+1}),\quad (t_1,\dots,t_{n+1})\in T^{n+1};\end{equation}
for each $(l_1,\dots,l_i)\in M$, $l_1+\dots,l_i+i=n$, and each $(t_1,\dots,t_n)\in T^{(l_1,\dots,l_i)}$,
\begin{equation}\label{dst} (X^0(f)g^{(n)})(t_1,\dots,t_n)=b^{(l_1)}(t_1)f(t_1)g^{(n)}(t_1,\dots,t_n);\end{equation}
and for each $(l_1,\dots,l_i)\in M$, $l_1+\dots,l_i+i=n-1$, and each $(t_1,\dots,t_{n-1})\in T^{(l_1,\dots,l_i)}$,
\begin{align}
(X^-(f)g^{(n)})(t_1,\dots,t_{n-1})&= \int_T \sigma(dt)\,f(t)g^{(n)}(t,t_1,\dots,t_{n-1})\notag\\
&\quad+ a^{(l_1+1)}(t_1)f(t_1)g^{(n)}(t_1,t_1,t_2,\dots,t_{n-1}) \label{trf}
\end{align}
(the second addend on the right hand side of \eqref{trf} being equal to zero for $n=1$). Additionally,
$X^+(f)\Omega=f$, $X^0(f)\Omega=0$, $X^-(f)\Omega=0$.

\end{theorem}

\begin{remark} For the reader's convenience, let us quickly summarize the constructed spaces and the established unitary isomorphisms. We first have the following commutative diagram:    

$$
\begin{diagram}
L^2(\tau)=\mathbb R\oplus\bigoplus_{(l_1,\dots,l_i)\in M}\mathbf H_{(l_1,\dots,l_i)} & &\rTo^I& & 
\mathcal F(\mathcal H)= \mathbb R\oplus\bigoplus_{(l_1,\dots,l_i)\in M}\mathcal H_{(l_1,\dots,l_i)} 
\\
 & \luTo_J &  & \ruTo_K & \\
 & & 
 \mathbb F=\mathbb R\oplus\bigoplus_{(l_1,\dots,l_i)\in M} \mathbb {H}_{(l_1,\dots,l_i)} 
 & & 
\end{diagram}
$$
Here, the spaces $\mathbb H_{(l_1,\dots,l_i)}$  are defined by \eqref{rtsryhgy}, the isomorphism $I$ is established in Proposition~\ref{yuf}, $K$ is given through  \eqref{yurftydeff},  $J$ is given by \eqref{vfytfytd}, the spaces $\mathcal H_{(l_1,\dots,l_i)}$ and $\mathbf H_{(l_1,\dots,l_i)}$ are the images of $\mathbb H_{(l_1,\dots,l_i)}$ under $K$ and $J$, respectively. Furthermore, we have realized 
each space 
$$ \mathbb F^{(n)}=\bigoplus_{(l_1,\dots,l_i)\in M,\ l_1+\dots+l_i+i=n} \mathbb {H}_{(l_1,\dots,l_i)},\quad n\in\mathbb N, $$
as $L^2(T^n,\gamma_n)$, and derived the following commutative diagram:
$$
\begin{diagram}
L^2(\tau)=\bigoplus_{n=0}^\infty\mathbf{OP}^{(n)} & &\rTo^I& & 
\mathcal F(\mathcal H) =\mathbb R\oplus\bigoplus_{n=1}^\infty \mathcal H^{(n)}
\\
 & \luTo_J &  & \ruTo_K & \\
 & & 
 \mathbb F=\mathbb R\oplus\bigoplus_{n=1}^\infty L^2(T^n,\gamma_n)
 & & 
\end{diagram}$$
 where 
 $$  \mathcal H^{(n)}:= \bigoplus_{(l_1,\dots,l_i)\in M,\ l_1+\dots+l_i+i=n} \mathcal {H}_{(l_1,\dots,l_i)},\quad n\in\mathbb N. $$
Formula \eqref{dyrrrs} gives the action of $J$ in terms of the latter diagram, while formulas \eqref{ydse}--\eqref{trf} give the action of $X(f)$ in $\mathbb F$. 
\end{remark}

\section{The free Meixner class}\label{gyder5asr}

As we saw in Theorem~\ref{ftty}, the free Gauss--Poisson  processes have the property that, for each $f^{(n)}\in\D^{(n)}$, the orthogonal polynomial $P(f^{(n)})$ is a continuous  polynomial. We will now search for all the free processes as in Section~\ref{serts}  for which this property remains true. So, as in Section~\ref{serts}, we fix a free process $(X(f))_{f\in\D}$ --- a family of bounded linear operators in  the free extended Fock  space $\mathbb F$.

\begin{theorem}\label{tyfxdzxsqx} The following statements are equivalent:

i) For each $f^{(n)}\in \D^{(n)}$, $ P(f^{(n)})\in\mathbf{C}{\mathbf P}$.

ii) For each $f\in\D$, $X(f)$ maps $\Ffin$ into itself.

iii) There exist $\lambda$ and $\eta$ from $C(T)$, $\eta(t)\ge0$ for all $t\in T$, such that
\begin{align*}
b^{(l)}(t)&=\lambda(t),\quad t\in T,\ l\in\mathbb N_0,\\
a^{(l)}(t)&=\eta(t),\quad t\in T,\ l\in\mathbb N.
\end{align*}
In this case, for each $f\in\D$ and $g^{(n)}\in\D^{(n)}$, $n\in\mathbb N$,
\begin{align}
&(X^+(f)g^{(n)})(t_1,\dots,t_{n+1})=f(t_1)g(t_2,\dots,t_{n+1}),\quad (t_1,\dots,t_{n+1})\in T^{n+1},
\label{eserssres}\\
&(X^0(f)g^{(n)})(t_1,\dots,t_n)=\lambda(t_1)f(t_1)g^{(n)}(t_1,t_2,\dots,t_n),\quad (t_1,\dots,t_n)\in T^n,\label{wewewe}\\
& (X^-(f)g^{(n)})(t_1,\dots,t_{n-1})
=
\int_T \sigma(dt)\,f(t)g^{(n)}(t,t_1,\dots,t_{n-1})\notag\\
&\qquad+ \eta(t_1)f(t_1)g^{(n)}(t_1,t_1,t_2,\dots,t_{n-1}),\quad (t_1,\dots,t_{n-1})\in T^{(n-1)}\label{yfy}
\end{align}
(the second addend on the right hand side of \eqref{yfy} being equal to zero for $n=1$).
\end{theorem}

\begin{proof}
Assume that i) holds. Hence, for any  $n\in\mathbb N$, there exist linear operators $U_{i,n}:\D^{(n)}\to\D^{(i)}$, $i=0,1,\dots,n$, such that
\begin{equation}\label{ydyt}  P(f^{(n)})=\sum_{i=0}^n\la U_{i,n}f^{(n)}, \omega^{\otimes i}\ra,\quad f^{(n)}\in\D^{(n)}.\end{equation}
Applying the orthogonal projection of $L^2(\tau)$ onto $\mathbf O{\mathbf P}^{(n)}$ to both right and left hand sides of \eqref{ydyt}, we get $P(f^{(n)})= P(U_{n,n}f^{(n)})$. Hence, $U_{n,n}$ is the identity operator, so that \eqref
{ydyt} becomes
\begin{equation}\label{ydftdftyt}  P(f^{(n)})=\la f^{(n)}, \omega^{\otimes n}\ra+\sum_{i=0}^{n-1}\la U_{i,n}f^{(n)}, \omega^{\otimes i}\ra,\quad f^{(n)}\in\D^{(n)}.\end{equation}
From here it follows that, for any  $n\in\mathbb N$, there exist linear operators $V_{i,n}:\D^{(n)}\to\D^{(i)}$, $i=0,1,\dots,n-1$, such that
\begin{equation}\label{ytdrtds}
\la f^{(n)},\omega^{\otimes n}\ra= P(f^{(n)})+\sum_{i=0}^{n-1}P( V_{i,n}f^{(n)}),\quad f^{(n)}\in\D^{(n)}.
\end{equation}
Indeed, for $n=1$, \eqref{ytdrtds} clearly holds. Assume that \eqref{ytdrtds} holds for all $n=1,\dots,N$, $N\in\mathbb N$. Then, by \eqref{ydftdftyt} and by \eqref{ytdrtds} for $n\le N$, we have,
for each $f^{(N+1)}\in\D^{(N+1)}$,
\begin{align*}
\la f^{(N+1)},\omega ^{\otimes(N+1)}\ra&= P(f^{(N+1)})-\sum_{i=0}^N\la U_{i,N+1}f^{(N+1)},\omega^{\otimes i}\ra\\
&=P(f^{(N+1)})-\sum_{i=0}^N\bigg(
 P(U_{i,N+1}f^{(N+1)})+\sum_{j=0}^{i-1} P(V_{j,i}U_{i,N+1}f^{(N+1)})
\bigg),
\end{align*}
from where \eqref{ytdrtds} holds for $n=N+1$.

Now, for each $f\in\D$ and $g^{(n)}\in\D^{(n)}$, by \eqref{ydftdftyt}  and \eqref{ytdrtds},
\begin{align}
&\la f,\omega\ra P(g^{(n)})=\la f\otimes g^{(n)},\omega^{(n+1)}\ra+\sum_{i=1}^{n}\la
f\otimes (U_{i-1,n}\,g^{(n)}),\omega^{\otimes i}\ra\notag\\
&\qquad= P(f\otimes g^{(n)})+\sum_{j=0}^n  P(V_{j,n+1}(f\otimes g^{(n)}))\notag\\
&\qquad\quad+\sum_{i=1}^{n}\bigg(
 P(f\otimes (U_{i-1,n}\,g^{(n)}))+\sum_{k=0}^{i-1} P(V_{k,i}(f\otimes (U_{i-1,n}\,g^{(n)}))
\bigg)\notag\\
&\qquad =P(f\otimes g^{(n)})+\sum_{j=0}^n  P (Z_{j,n+1}(f, g^{(n)})),\label{dtrdxstr}
\end{align}
where $Z_{j,n+1}(f, g^{(n)})\in \D^{(j)}$.
Thus, by \eqref{dtrdxstr}, $\la f,\omega\ra P(g^{(n)})\in\mathbf C{\mathbf P}^{(n+1)}$, and so ii) holds. (Note that, in view of symmetricity, $Z_{j,n+1}(f,g^{(n)})=0$ for $j\le n-2$.)

Let us now prove that ii) implies iii). For each $t\in T$, denote $\lambda(t):=b^{(0)}(t)$ and $\eta(t):=a^{(1)}(t)$.

Fix any open set $O\in\mathcal B_0(T)$. Let $f,g\in\D$ be such that $f(t)=g(t)=1$ for all $t\in O$. Then, by
\eqref{dst}, for each $t\in O$,
\begin{equation}\label{hjy}
(X^0(f)g)(t)=\lambda(t).
\end{equation}
By ii), \eqref{hjy} implies that $\lambda(t)$ continuously depends on $t\in O$. Hence, $\lambda\in C(T)$.

Next, let a set $O$ and functions $f$, $g$ be as above, and assume additionally that $f\ge0$ and $g\ge0$ on $T$. Further, choose any $h\in\D$ such that $h\ge0$ on $T$, with  $h(t)=0$ for all $t\in O$, and $fh\not\equiv0$.  Choose $\varepsilon<0$ such that
$$ \int_T\sigma(dt)f(t)(g(t)+\varepsilon h(t))=0. $$
Set $g^{(2)}(t_1,t_2):=(g(t_1)+\varepsilon h(t_1))g(t_2)$, $(t_1,t_2)\in T^2$. Then, by \eqref{trf}, for each $t\in O$,
$$ (X^-(f)g^{(2)})(t)=\eta(t),$$ which implies
that $\eta$ is continuous on $O$. Hence, $\eta\in C(T)$.

Next, fix any $t\in T$, and let $f\in\D$ and $g^{(n)}\in\D^{(n)}$, $n\ge2$, be such that
$f(t)=1$ and $g^{(n)}(t,t,\dots,t)=1$.
 By \eqref{dst},
for any $(t_1,\dots,t_{n-1})\in T^{n-1}$ such that
$t\ne t_1$, $t_1\ne t_2$, \dots, $t_{n-2}\ne t_{n-1}$, we have
\begin{equation}\label{rdrt}
(X^0(f)g^{(n)})(t,t_1,t_2\dots,t_{n-1}) =
\lambda(t)g^{(n)}(t,t_1,t_2,\dots,t_{n-1}),\end{equation}
whereas
\begin{equation}\label{ydjkd}
(X^0(f)g^{(n)})(t,t,\dots,t)=b^{(n-1)}(t).
\end{equation}
By ii),
$$\lim_{(t_1,t_2,\dots,t_{n-1})\to(t,t,\dots,t)}(X^0(f)g^{(n)})(t,t_1,t_2,\dots,t_{n-1})=(X^0(f)g^{(n)})(t,t,\dots,t).$$
Hence, by \eqref{rdrt} and \eqref{ydjkd}, $b^{(n-1)}(t)=\lambda(t)$. Thus, for all $t\in T$ and all $n\in\mathbb N_0$, $b^{(n)}(t)=\lambda(t)$.

Completely analogously, we then also deduce from \eqref{trf} that, for all $t\in T$ and all $n\in\mathbb N$, $a^{(n)}(t)=\eta(t)$.
Formulas \eqref{wewewe} and \eqref{yfy} now follow from \eqref{dst} and \eqref{trf}, respectively.
Thus, iii) holds.

Finally, we prove that iii) implies i).
Analogously to \eqref{tyfsqx}, we now have, for any $f_1,\dots,f_n\in\D$, $n\ge2$:
\begin{align*}
& P(f_1\otimes\dots\otimes f_n)=\la f_1,\omega\ra P(f_2\otimes\dots\otimes f_n)-
 P((\lambda f_1 f_2)\otimes f_3\otimes\dots\otimes f_n)\\
&\quad-\int_T\sigma(dt)f_1(t)f_2(t) P(f_3\otimes\dots\otimes f_n)- P((\eta f_1f_2f_3)\otimes f_4\otimes\dots\otimes f_n)
\end{align*}
(compare with \eqref{tyfsqx}).
From here we conclude statement i) by an easy generalization of the proof of Theorem~\ref{ftty}.
\end{proof}

\begin{remark} \label{hgvufv} As easily seen by approximation, formulas \eqref{eserssres}--\eqref{yfy} remain true for any $f\in B_0(T)$ and $g^{(n)}\in\mathbb F^{(n)}=L^2(T^n,\gamma_n)$.
\end{remark}

The set of all free processes as in Theorem~\ref{tyfxdzxsqx}, iii) will be called the Meixner class of free processes. We note that, if, for $t\in T$, $\eta(t)=0$, then the measure $\mu(t,\cdot)$ is concentrated at one point, namely $\lambda(t)$. Hence, $g^{(0)}(t)=1$ and $g^{(l)}(t)=0$ for all $l\in\mathbb N$ (see \eqref{tdtdtdftt} and \eqref{sresydrtdstrs}).
In particular, if $\eta(t)=0$ for all $t\in T$,
the measure $\gamma_n$ becomes $\sigma^{\otimes n}$ (see \eqref{11111}). Thus, $\mathbb F=\mathcal F(\mathcal H)$ and $X(f)=x(f)$, $f\in\D$, where $(x(f))_{f\in\D}$ is the free process as in Section~\ref{cfrd}, which  corresponds to the function $\lambda\in C(T)$.

If, however, $\eta(t)>0$, then $\mu(t,\cdot)$ has an infinite support. Recall that $\mu(t,\cdot)$ is the measure of orthogonality of monic polynomials $(p^{(n)}(t,\cdot))_{n=0}^\infty$ satisfying
\begin{equation}\label{yutto8} sp^{(n)}(t,s)=p^{(n+1)}(t,s)+\lambda(t)p^{(n)}(t,s)+\eta(t) p^{(n-1)}(t,s),\quad n\in\mathbb N_{0},\end{equation}
where $p^{(-1)}(t,s):=0$. Hence,  $\mu(t,\cdot)$ is Wigner's semicircle law with mean $\lambda(t)$ and variance $\eta(t)$:
$$ \mu(t,ds)=\chi_{[-2\sqrt{\eta(t)}+\lambda(t),\,2\,\sqrt{\eta(t)}+\lambda(t)]}(s)\,(4\pi\eta(t))^{-1}\,\sqrt{4\eta(t)-(s-\lambda(t))^2}\,ds$$
(compare with \cite{SY} and \cite{BB}).
By \eqref{tdtdtdftt} and \eqref{yutto8}, we have: 
\begin{equation}\label{tydetys} g^{(l)}(t)=\eta^l(t),\quad l\in\mathbb N_0.
\end{equation}
Substituting \eqref{tydetys} into \eqref{iutfrsdr}, we get the explicit form of the inner product in the free extended Fock space $\mathbb F$. 

Assume that, for some $\Delta\in\mathcal B_0(T)$, the functions $\lambda(\cdot)$ and $\eta(\cdot)$ are constant on $\Delta$, i.e., $\lambda(t)=\lambda$, $\eta(t)=\eta$ for all $t\in \Delta$, where $\lambda\in\mathbb R$ and $\eta\ge0$. Then, by \eqref{eserssres}--\eqref{yfy} (see Remark~\ref{hgvufv}), we have:
\begin{equation}\label{xzzdz}X(\Delta)\chi_\Delta^{\otimes n}=
\chi_{\Delta}^{\otimes(n+1)}+[n]_0 \lambda\chi_\Delta^{\otimes n}+([n]_0\sigma(\Delta)+[n]_0[n-1]_0\eta)\chi_\Delta^{\otimes(n-1)},\quad n\in\mathbb N_0,\end{equation}
where $\chi_\Delta^{\otimes0}:=\Omega$. Denote
$ P(\chi_\Delta^{\otimes n}):=J\chi_\Delta^{\otimes n}$. Then, by \eqref{xzzdz},
$P(\chi_\Delta^{\otimes n})=q^{(n)}(X(\Delta))$,
where $(q^{(n)})_{n=0}^\infty$ is the system of monic polynomials on $\mathbb R$ satisfying the recursive relation \eqref{yted6}. 
By Favard's theorem, $(q^{(n)})_{n=0}^\infty$ is a system of polynomials which are orthogonal with respect to some probability measure $\rho_{\lambda,\eta,\sigma(\Delta)}$. For an explicit form of this measure, we refer to e.g.\ \cite{SY}.

\begin{corollary}\label{tdrssd}
Let $(X(f))_{f\in B_0(T)}$ be as in Theorem~\ref{tyfxdzxsqx} iii). Then, for each $\Delta\in\mathcal B_0(T)$, there exists $r=r(\Delta)>0$ such that, for each $f\in B_0(T)_{\mathbb C}$ satisfying
$$ |f(t)|<r\chi_{\Delta}(t)\quad\text{for all }t\in T, $$ we have
$$C(f)=\int_T\sigma(dt)\,2f^2(t)\left(1-\lambda(t)f(t)+\sqrt{(1-\lambda(t)f(t))^2-4f^2(t)\eta(t)}\right)^{-1}\,.$$
\end{corollary}

\begin{proof} The result directly follows from Proposition~\ref{hdrtdx} and the following formula which holds for $z\in\mathbb C$ from a neighborhood of zero:
$$ \int_{\mathbb R}\mu(t,ds)\frac1{1-sz}=2\left(
1-\lambda(t)z+\sqrt{(1-\lambda(t)z)^2-4z^2\eta(t)}
\right)^{-1},$$
see \cite{a1,SY}.
\end{proof}

Recall that $\Ffin$ is a dense subset of $\mathbb F$. Analogously to Section~\ref{cfrd}, we can therefore interpret smeared, Wick ordered products of operators $\di_t^\dag$
and $\di_t$ as operators in $\mathbb F$. 

\begin{corollary}\label{gfdctdx}
Let $(X(f))_{f\in B_0(T)}$ be as in Theorem~\ref{tyfxdzxsqx} iii). Then, using the same notations as in Section~\ref{cfrd}, we may represent the action of each $X(f)$ in $\mathbb F$ as follows:
$$ X(f)=\int_T\sigma(dt)f(t)\omega(t),$$
where
$$ \omega(t)=\di_t^\dag+\lambda(t)\di_t^\dag\di_t+\di_t+\eta(t)\di_t^\dag\di_t\di_t.$$
\end{corollary}

\begin{proof}
The statement directly follows from \eqref{eserssres}--\eqref{yfy} if we note that, for each $g^{(n)}\in\D^{(n)}$,
$$\bigg(\int_T\sigma(dt)f(t)\eta(t)\di_t^\dag\di_t\di_t g^{(n)} \bigg)(t_1,\dots,t_{n-1})=\eta(t_1)f(t_1)g^{(n)}(t_1,t_1,t_2,\dots,t_{n-1}) .$$
\end{proof}

\begin{center}
{\bf Acknowledgements}\end{center}

 We would  like to thank the referee for a careful reading of the
manuscript and making very useful comments and suggestions. The authors acknowledge the financial support of the SFB~701 ``Spectral structures and topological methods in mathematics'', Bielefeld University. MB was partially supported by the  KBN grant no.\ 1P03A 01330. EL was partially supported by  the PTDC/MAT/67965/2006 grant, University of Madeira.

\end{document}